\documentclass[journal]{IEEEtran}

\usepackage [english]{babel}
\usepackage [autostyle, english = american]{csquotes}
\MakeOuterQuote{"}

\IEEEoverridecommandlockouts

\usepackage[pdftex]{graphicx}
\usepackage{epstopdf}
\usepackage{amsmath}
\usepackage{amsfonts}
\usepackage{balance}

\usepackage{amssymb}
\usepackage{bbm}

\usepackage{array}
\usepackage{multirow}
\usepackage{color}
\usepackage{tablefootnote}

\usepackage[prependcaption,colorinlistoftodos]{todonotes}

\addto\captionsenglish{}

\newtheorem{defn}{Definition}
\newtheorem{thm}{Theorem}
\newtheorem{proof}{Proof}

\newtheorem{rmk}{Remark}
\newtheorem{proposition}{Proposition}

\newtheorem{assumption}{Assumption}

\newcommand{\jw}{\textcolor{black}}
\newcommand{\jdw}{\textcolor{black}}
\newcommand{\icl}{\textcolor{black}}
\newcommand{\jd}{\textcolor{black}}
\newcommand{\il}{\textcolor{black}}
\newcommand{\jwjw}{\textcolor{black}}
\newcommand{\ill}{\textcolor{black}}
\newcommand{\jwjwjw}{\textcolor{black}}
\newcommand{\jdjw}{\textcolor{black}}
\newcommand{\ic}{\textcolor{black}}
\newcommand{\jww}{\textcolor{black}}
\newcommand{\icc}{\textcolor{black}}
\newcommand{\jwww}{\textcolor{black}}
\newcommand{\jwwww}{\textcolor{black}}
\newcommand{\iclc}{\textcolor{black}}
\newcommand{\jjww}{\textcolor{black}}
\newcommand{\jjdw}{\textcolor{black}}
\newcommand{\jjjw}{\textcolor{black}}

\newcommand{\iclcl}{\textcolor{black}}
\newcommand{\ilc}{\textcolor{black}}
\newcommand{\ilcl}{\textcolor{black}}

\definecolor{klcolor}{rgb}{0.0,0.5,0}
\newcommand{\kl}{\textcolor{black}}

\newcommand{\li}{\textcolor{black}}
\newcommand{\lic}{\textcolor{black}}
\newcommand{\licc}{\textcolor{black}}
\newcommand{\lif}{\textcolor{black}}

\newcommand{\kkl}{\textcolor{black}}

\newcommand{\jwrev}{\textcolor{black}}
\newcommand{\jdwrev}{\textcolor{black}}
\newcommand{\jdwr}{\textcolor{black}}
\newcommand{\jwr}{\textcolor{black}}

\newcommand{\jjwr}{\textcolor{black}}

\newcommand{\jwwj}{\textcolor{black}}

\begin{document}

\title{A scalable control design for grid-forming inverters in microgrids}

\author{Jeremy D.~Watson, Yemi Ojo, Khaled Laib, and Ioannis Lestas
\thanks{J.D.~Watson, Yemi Ojo, Khaled Laib, and Ioannis Lestas are with the Department of Engineering, University of Cambridge, 
\ilc{Cambridge CB21PZ,} the United Kingdom. Emails: {jdw69, yo259, kl507, icl20} @cam.ac.\iclcl{uk.}}}
\maketitle

\begin{abstract}

Microgrids are increasingly recognized as a key technology for the integration of \jdw{distributed energy resources} into the power network, allowing local clusters of load and distributed energy resources to operate autonomously.  However, microgrid operation brings new challenges, especially in islanded operation as frequency and voltage control \jdwr{are} no longer provided by large rotating machines. Instead, the power converters in the microgrid must coordinate to regulate the frequency and voltage and ensure stability. We consider the problem of designing controllers to achieve these objectives. \jdw{Using passivity theory
to derive decentralized stability conditions for the microgrid}, we propose a control design method for grid-forming inverters.
\ic{For the analysis we use higher-order models for the inverters and also advanced dynamic models for the lines with an arbitrarily large number of states}. \jdw{By satisfying the \ic{decentralized condition formulated}}, plug-and-play operation can be achieved with guaranteed stability, \ic{and performance can also be improved by incorporating this condition as a constraint in corresponding optimization problems formulated.} In addition, our control design 
can improve the power sharing \ic{properties} 
of the microgrid compared to previous non-droop approaches. Finally, realistic simulations confirm that the controller design improves the stability and performance of the power network.

\end{abstract}



\IEEEpeerreviewmaketitle

\section{Introduction}

Microgrids are a key technology for the integration of \jdw{distributed energy resources (DER)} into the power network. These distributed energy resources are connected to the network via power converters. The microgrid concept allows local clusters of loads, distributed generation, and energy storage systems to operate together in smaller electrical networks, either in grid-connected mode or autonomously. The potential benefits are considerable, including reduced capital costs, increased efficiency and reliability. However, microgrid operation brings new challenges of its own, since power converters have contrasting dynamical behavior to rotating machines with large inertia~\cite{milano2018}. In islanded operation, frequency and voltage control are no longer provided by large rotating machines. Instead, the power converters in the microgrid must \jdw{regulate the frequency and voltage themselves, while ensuring stability}. Furthermore, they must \jdw{share the power requirement} between the converters to avoid exceeding ratings. 

Traditionally, droop control has been used for \jdw{inverters in} microgrids, usually in terms of active power-frequency droop and reactive power-voltage droop~\cite{chandorkar1993}. These take their inspiration from \jwrev{similar} mechanisms built into the governors and exciters of synchronous generators. Droop controllers are practical, simple, and operate entirely on local measurements. However, they suffer from several issues, including frequency and voltage deviations depending on the load \ilc{and} poor transient performance \jdw{in some} situations~\cite{yao2011}~\cite{watson2019}. 
\ilc{Analytical studies also often incorporate additional assumptions, such as simplified inverter models or static/lossless lines, to guarantee stability
in a decentralized way in general network topologies.} Another possible droop approach for stabilizing inverter-based power systems is angle droop control~\cite{majumder2009}, and the closely related virtual impedance technique~\cite{sun2017}. In our previous work~\cite{watson2019}, these schemes were compared and the difficulty of guaranteeing 
\iclc{stability 
in these schemes via decentralized passivity-based conditions} was shown.

Non-droop approaches typically require a Power Measurement System (PMS) which solves an optimal power flow (OPF) problem for the entire microgrid at specified intervals and communicates appropriate setpoints to each inverter. While this requires additional communication between the PMS and each inverter, better transient performance may be achieved along with optimal steady-state power allocations. The accuracy of the load sharing is dependent on how often the optimal power flow is solved~\cite{sadabadi2017}, and the grid frequency is fixed by internal oscillators which are periodically synchronized. Proposed approaches include $H_\infty$ controllers~\cite{derakhshan2020}, two-degrees-of-freedom control~\cite{babazadeh2013}, robust servomechanism controllers~\cite{etemadi2014}, \jw{and LMI-based state feedback~\cite{tucci2020}. However, \jd{some} of these approaches require retuning to integrate new devices and \jwrev{many} do not take the dynamical nature of the distribution lines into account.} In~\cite{sadabadi2017} polytopic uncertainty was used to achieve plug-and-play capability without this requirement for retuning; however, an assumption that certain coupling terms between the inverters can be neglected is required, which in practice may not always be satisfied. \jwrev{In~\cite{strehle2019}, a port-Hamiltonian formulation was used to achieve plug-and-play stabilization of the microgrid. One general drawback of \jdwrev{many of \lic{the} non-droop approaches}} is that 
during the period between the load change and the new OPF broadcast, proportional power sharing cannot be guaranteed~\cite{sadabadi2017}, and these schemes are usually unable to successfully operate in the case of communication failure. 

\jdw{In this paper, we aim to combine the advantages of both droop and non-droop approaches by \jdwr{\li{means of a} control design} which, similarly to the non-droop controllers, is synchronized to the other inverters in the network via a PMS. \jd{However, by designing a virtual impedance that exemplifies angle droop-like behavior~\cite{sun2017}, we also allow the power requirements to be shared in the time between setpoint broadcasts}. \jwrev{Compared to previous work~\cite{he2018} using virtual impedance and a PMS, we design novel controllers for which we are able to guarantee stability in a decentralized way. }
\jdwrev{We} \jdw{use passivity theory \jwrev{for this stability guarantee.}}} In order to apply passivity approaches to AC power networks, the choice of an appropriate frame of reference for representing the network dynamics requires consideration. \jdw{Many approaches \ilc{for stability analysis}} 
have focused on models formulated in a local reference frame. \jdw{However, \ilc{for the network to be passive} this typically requires the assumption of lossless static lines.} \ilc{An alternative} 
approach is to use the Park and Clarke transformations in a common reference frame~\cite{spanias2018}, since in this frame of reference transmission lines retain their natural passivity properties, without having to resort to a lossless assumption or omit their dynamic \jdwr{behavior}. This frame of reference is often used in the literature for grid-connected voltage-source converters (VSCs), e.g.~\cite{harnefors2016}, \il{and its significance in inverter-based grids was \ic{also} pointed out in~\cite{watson2019}}. \jwrev{In terms of microgrids with grid-forming inverters, the key contribution of~\cite{nahata2019} showed that the same passivity property leads to network stability. Compared to this and other work in the literature, our virtual impedance approach \jjwr{allows to establish strict passivity while providing power sharing. Furthermore, we \li{exploit} the LMI formulation for passivity \li{to form\jwwj{ulate} optimization problems} which can lead to static as well as high-order dynamical controllers that satisfy additional performance criteria.}}

\subsection*{\kl{Paper contributions}}
\jdw{This paper 
\icl{considers} the problem of \jdwr{control design} for grid-forming inverters which guarantees stability and plug-and-play operation via passivity in a common reference frame, while still having appropriate properties that facilitate power-sharing and voltage regulation. }
The main contributions of this work may be summarized as follows:
\begin{enumerate}
    \item \jw{We propose a 
    \il{control design for grid-forming inverters} which \iclc{ensures} plug-and-play \ilc{capability by} 
    \il{satisfying a} decentralized passivity criterion. 
     \il{In particular, the proposed \jdwr{controllers provide} stability guarantees for the microgrid}
        \icc{
        using a general dynamic model for the \ilc{lines and a higher-order model for the inverters.}}}
    \item \jdwr{The use of passivity as an LMI constraint allows us to formulate optimization problems with various objectives or other constraints to improve performance. This can lead to various controllers that preserve the passivity property, such as \li{a dynamic controller that satisfies 
        mixed $H_\infty$\jwwj{/}passivity conditions.}} 
    \item \jw{Our control strategy, via its droop-like behavior, allows contributions from multiple inverters to regulate disturbances even before new setpoints are broadcast.}
    \item \jdw{A case study using detailed simulations is used to verify the performance of the proposed \jdwr{controllers}, demonstrating the \icl{very good} 
        performance of an inverter-based microgrid with the proposed \jdwr{controllers}.}
\end{enumerate}

\subsubsection*{\kl{Paper outline}}
The paper is organized as follows: in Section~\ref{prelim} we present preliminaries, introduce the models used, and discuss the passivity condition for the stabilization of the network. \jdw{In Section~\ref{design} we present our control design.} 
Simulations presented in Section~\ref{casestudy} show how the dynamic response of the system is improved by our passivity-based design. Finally, Section~\ref{conclusion} summarizes our conclusions.

\section{Preliminaries and models}\label{prelim}

\subsection{Notation and definitions}\label{notdef}

A power system is comprised of buses (generators, loads, etc.) and lines, which we represent as a graph $(N, E)$ where $N$ is the set of buses $\{1,2, \dots |N|\}$ and $E \subset N \times N$ is the set of edges (lines). \jd{A direction is assigned to each edge which can be arbitrarily chosen. }
%

In this paper, we use the well-known Park and Clarke transformations to transform a balanced AC signal $x_{abc}$ into its $DQ$ components in a common reference frame with \jdw{constant} synchronous frequency $\omega_{\jdjw{s}}$. The transformation $T(\omega_{\jdjw{s}} t)$ we adopt is the power invariant \ilc{form,}
which aligns the D-axis with the inverter terminal voltage.
\begin{gather*}
	T(\omega_{\jdjw{s}} t)=\sqrt{\frac{2}{3}}
	\begin{bmatrix}
	\cos(\omega_{\jdjw{s}} t) &\cos(\omega_{\jdjw{s}} t-\frac{2\pi}{3}) &\cos(\omega_{\jdjw{s}} t+\frac{2\pi}{3}) \\
	\sin(\omega_{\jdjw{s}} t) &\sin(\omega_{\jdjw{s}} t-\frac{2\pi}{3}) &\sin(\omega_{\jdjw{s}} t+\frac{2\pi}{3})
 \end{bmatrix}\\
    \ilc{x_{DQ}=T(\omega_{\jdjw{s}} t)x_{abc}}
\end{gather*}
\jd{The voltage at a bus $j \in N$ is denoted by $v^j_{DQ}$, the summation of the currents flowing out of that bus into one or more edges (i.e. \jdwr{what is often called} the ``current injection'' \jdwr{of} that bus \jdwr{to the network}) is denoted by $i^j_{DQ}$.}
We will use \ilc{$V_{DQ} = [v^j_{DQ}]_{j\in N}$, $I_{DQ} = [i^j_{DQ}]_{j\in N}$ to denote the vectors of all DQ bus voltages and current injections, respectively.}
$I_n$ denotes the $n\times n$ identity matrix, \jwjwjw{$J$ the matrix $\begin{bmatrix}0 & 1\\ -1 & 0\end{bmatrix}$,} $\otimes$ denotes the Kronecker product, \jwwj{and $\Re[z]$ denotes the real part of \lif{a complex} number $z$}. Throughout the paper, the \ic{equilibrium} value of a variable $x$ is denoted by an asterisk \il{superscript $x^*$,} \jwwww{\ilc{the} eigenvalues of a matrix $A$ are denoted by $\lambda_i(A)$ and the maximum singular value by $\bar{\sigma}(A)$.} \jjwr{We use $\text{diag}(x)$ to denote the diagonal matrix with the elements of the vector $x$ on the diagonal, and $\textbf{bdiag}$ to denote the \li{analogous operation} 
for block diagonal matrices.}

\subsection{Line dynamics}\label{nw}

\begin{figure}[ht]
    \centering
    \includegraphics[width=0.5\textwidth]{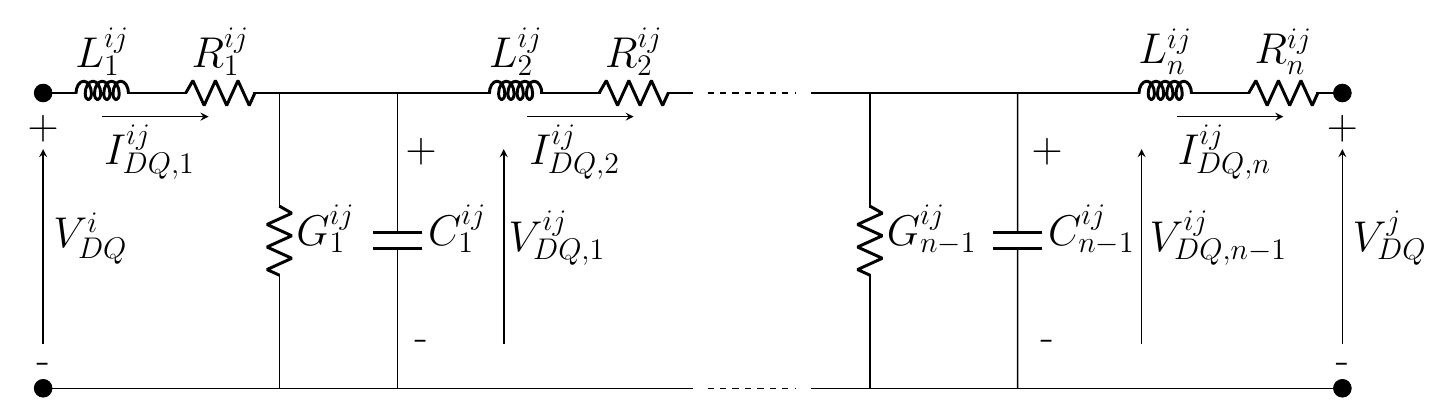}
    \caption{\jwjwjw{Transmission line \iclcl{model.}}}
    \label{fig:linemodel}
\end{figure}

We model the transmission/distribution lines as RLC components and assume that the network is three-phase balanced in all lines and at all buses. The dynamical equations in the $DQ$ frame are easily derived by applying the Park-Clarke transformation to the fundamental equations of the \iclc{individual components.} 
\jwjw{We write these equations}
\ill{in a general form, where lines are represented with a distributed parameter model with \iclc{an} arbitrarily large number of states and with the line parameters allowed to vary along the lines, \jwjwjw{as illustrated in Fig.~\ref{fig:linemodel}}. It should be noted that the simpler lumped model for transmission lines is a special case of this model and it is also known that this tends to the telegrapher equations as the number of elementary sections tends to infinity\footnote{The lumped parameter model is generally sufficient for short lines, however, we include the more detailed model for completeness. The distributed model may be relevant for analysis at fast timescales where the wavelength is short, or in cases where parameters vary along the lines.}.} 

\iclcl{In particular, each line $(i, j) \in E$ has $n$ subsections, where $n$ can be an arbitrarily large number \ilc{(the stability conditions that will be presented are independent of $n$).}
Each subsection $k$ has resistance/inductance/conductance/capacitance $R^{ij}_k, L^{ij}_k, G^{ij}_k, C^{ij}_k$ respectively, as shown in Fig.~\ref{fig:linemodel},  and the line equations in DQ coordinates are:}
\begin{equation}
	{\small
\begin{matrix}
    L^{ij}_1\dot{I}^{ij}_{DQ,1}&=& -R^{ij}_1I^{ij}_{DQ,1} + \jwww{\omega_{s}}L^{ij}_1JI^{ij}_{DQ,1} + V^i_{DQ} - V^{ij}_{DQ,1} \\
    C^{ij}_1\dot{V}^{ij}_{DQ,1} & = &-G^{ij}_1V^{ij}_{DQ,1} + \omega_{s}C^{ij}_1JV^{ij}_{DQ,1} + I^{ij}_{DQ,1} - \ilc{I^{ij}_{DQ,2}}\\
    &\vdots&\\
    L^{ij}_n\dot{I}^{ij}_{DQ,n}& =&-R^{ij}_nI^{ij}_{DQ,n} + \omega_{s}L^{ij}_nJI^{ij}_{DQ,n} + V^{ij}_{DQ,n-1} - V^j_{DQ}\\[2ex]
\end{matrix}}\\[1ex]
\label{eq:linemodel}
\end{equation}
where the inductor current of the $k$th subsection of the line $(i,j) \in E$ is $I^{ij}_{DQ,k}$, the capacitor voltage of the $k$th capacitor of the same line is $V^{ij}_{DQ,k}$, and the terminal (end-point) voltages are $V^i_{DQ}$ and $V^j_{DQ}$ where $i$ and $j$ are the buses which are connected by the line $(i,j)$.

\jwrev{By defining various incidence matrices \jdwr{$H, H_B, H_C$ (detailed in  \kl{Appendix A}) composed of zeros and ones} to simplify the notation,}~\eqref{eq:linemodel} for all lines $(i, j) \in E$ can be represented in matrix form as follows:
\begin{equation}\label{netwdyn}
{\small \begin{matrix}
    [\mathcal{L} \otimes I_2]\dot{I}^\mathcal{N}_{dq} &=& -[\mathcal{R} \otimes I_2 + \omega_s\mathcal{L} \otimes J]I^\mathcal{N}_{dq} + [H_C^T \otimes I_2]V^C_{dq}\\
     &&\phantom{hi}+[H_B^T \otimes I_2]V_{dq} \\
    [\mathcal{C} \otimes I_2]\dot{V}^C_{dq} &=& -[\mathcal{G} \otimes I_2 + \omega_s\mathcal{C} \otimes J]V^C_{dq} + [H_C \otimes I_2]I^\mathcal{N}_{dq}\\
    I_{dq} &=& [H_B \otimes I_2]I^\mathcal{N}_{dq}
\end{matrix}}\\[2ex]
\end{equation}
{where $\omega_{\jdjw{s}}$ is the constant synchronous frequency, and \jdwr{$\mathcal{R}, \mathcal{L}, \mathcal{G}, \mathcal{C}$ are diagonal matrices of the \jwr{line} resistances, inductances, conductances, and capacitances respectively.}} \jwjwjw{As stated in Section~\ref{notdef},} $I_{DQ}$ and $V_{DQ}$ represent the vector of currents drawn from each bus, and the bus voltages respectively.

Note that at each bus $j$, the current injected into the network $i^j_{DQ}$ and the bus voltage $v^j_{DQ}$ are related by the bus dynamics, which in this study will include the inverter dynamics discussed in Section~\ref{gfc}.

\subsection{\jdwr{Bus} dynamics} \label{gfc}
  \kl{We consider now the bus dynamics. As mentioned earlier, there \jwr{are} two types of buses: generation and load.  The main focus of this paper is the generation buses where we \jdwr{design controllers} to ensure stability together with some performance requirements. \jjwr{We discuss allowable load buses and load dynamics in Section \ref{passdefn}, Remark \ref{loadrmk}.}} 
  The converter dynamics \kl{at generation buses} are allowed to be different at each bus $j$, however the bus index $j$ will be omitted for convenience in the presentation.
\begin{assumption}\label{invass}
We use an averaged inverter model with a stiff DC voltage. We also assume that harmonics can be ignored.
\end{assumption}

\begin{figure}[t]
    \centering
    \includegraphics[width=0.3\textwidth]{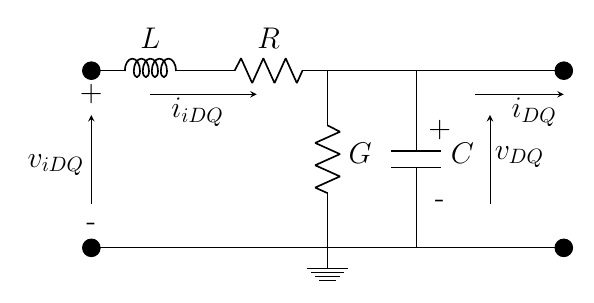}
    \caption{Inverter filter circuit diagram, showing the notation and current direction convention used in this paper. The inverter itself is connected on the left, with $v_{iDQ}$ representing the inverter voltage behind the filter.}
    \label{fig:iofilt}
\end{figure}

The inverter models in this paper are LC-filtered. We use $v_{iDQ}$ to refer to the inverter voltage which is the control input to the inductive filter, and $v_{DQ}$ is the capacitor voltage. Similarly, the current entering the LC filter \jdw{from the DC-side} is $i_{iDQ}$ and the current leaving the filter (injected into the network) is $i_{DQ}$. Fig.~\ref{fig:iofilt} \ilc{illustrates} 
the voltages and currents used in the model.

Using Assumption~\ref{invass}, the converter dynamics may be represented as a linear system of the form:
\begin{equation}\label{eq:convdyn}
\begin{matrix}
    \dot{x} &= \mathcal{A}x+\mathcal{B}_uu + \mathcal{B}_ww\\
    z&=\kl{\mathcal{C}_z}x+\mathcal{D}_uu+\mathcal{D}_ww
\end{matrix}
\end{equation}
where the states are $x = [i_{iDQ}^T \phantom{i} v_{DQ}^T]^T$, $u=v_{iDQ}$ is the controlled PWM voltage, and the input current from the network lines is \jdw{the} current $w = -i_{DQ}$ and the output is the voltage \jwrev{vector} $z = v_{DQ}$. The system matrices are given by:
\begin{equation}\label{eq:ssABCD}{\small
		\begin{aligned}
\mathcal{A} &= \begin{bmatrix}-\dfrac{R}{L} & \omega_{\jdjw{s}} & -\dfrac{1}{L} & 0\\ -\omega_{\jdjw{s}} & -\dfrac{R}{L}& 0 & -\dfrac{1}{L}\\ 1/C & 0 & -\dfrac{G}{C} & \omega_{\jdjw{s}}\\ 0 & \dfrac{1}{C} & -\omega_{\jdjw{s}} & -\dfrac{G}{C}  \end{bmatrix},
\mathcal{B}_u = \begin{bmatrix}\dfrac{1}{L} & 0\\0 & \dfrac{1}{L} \\ 0 & 0 \\ 0 & 0 \end{bmatrix}\\[2ex]
\mathcal{B}_w& = \begin{bmatrix}0 & 0\\0 & 0 \\ \dfrac{1}{C} & 0 \\ 0 & \dfrac{1}{C} \end{bmatrix},
\kl{\mathcal{C}_z} = \begin{bmatrix}0 & 0 & 1 & 0 \\ 0 & 0 & 0 & 1\end{bmatrix},
\mathcal{D}_u = \mathcal{D}_w = \begin{bmatrix}0 & 0\\0 & 0\end{bmatrix}
        \end{aligned}}\\[2ex]
\end{equation}

In order for each inverter to respond to broadcast signals from a PMS, \jdwr{the dynamics may be augmented with} an integrator of the form:
\begin{equation}\label{eq:int}
    \dot{\zeta}_{DQ} = v_{DQ}-v_{DQ}^{set}+Zi_{DQ}
\end{equation}
where $v_{DQ}^{set}$ is \jjwr{the voltage setpoint from the PMS} and $Z$ is a $2 \times 2$ matrix of coefficients which determine the \jjwr{terminal} voltage as a function of the input current. Previous work, e.g.~\cite{sadabadi2017}, used pure integral action (i.e $Z = 0_{2\times 2}$) to achieve this voltage regulation. In our \jdwr{full-state feedback} design we treat $Z$ as a virtual impedance with $Z = R_VI_2 - X_VJ$,
\jwjw{where $R_V > 0$ and $X_V > 0$ are the virtual resistance and reactance, respectively.} \icc{We denote the augmented dynamics as}:
\begin{equation}\label{eq:convdyn2}
\begin{matrix}
    \dot{\tilde{x}} &= \tilde{A}\tilde{x}+\tilde{B}_uu + \tilde{B}_ww\\
    z&=\tilde{C}\tilde{x}+\tilde{D}_uu+\tilde{D}_ww
\end{matrix}
\end{equation}
where $\tilde{x} = [i_{iDQ}^T \phantom{i} v_{DQ}^T \phantom{i} \zeta_{DQ}^T]^T$ \jwww{and the system matrices follow from \jww{\eqref{eq:ssABCD}-\eqref{eq:int}}.} Our approach using~\eqref{eq:int} has several advantages over previous approaches. Firstly, \jdwrev{power} sharing is achieved during unscheduled load changes. The \jwrev{terms} in $Z$ \jwjw{are a virtual impedance that allows the power requirement of load changes to be shared among the converters.} \jdwrev{The power sharing depends on the equivalent impedance at each converter~\cite{sun2017}. Supposing the line impedances to be comparatively small, \jww{the power allocation will depend mainly on the virtual impedances. Since these can be chosen by design, an appropriate allocation can thus be} achieved. \jwrev{In the case that the line impedances are significant, the equivalent impedance at each inverter determines its contribution and \jwr{a slight} inaccuracy will occur, as was analyzed in~\cite{he2018} for example.}} This can be adjusted to reach an acceptable compromise between power sharing and voltage fluctuations~\cite{he2018}. \jwrev{Secondly, our approach also allows} the design of a controller \icl{such that a local strict passivity condition is satisfied}  
without excessively high control gains. \jdw{With \icl{a pure integrator}, it is actually impossible to achieve strict passivity~\cite{nahata2019}.}

\jwrev{Note that any power sharing inaccuracy is temporary as no inaccuracy is introduced into the OPF by our proposed approach as long as the virtual impedance $Z$ at each inverter is known to the PMS.} The voltage setpoints $v_{DQ}^{set}$ broadcast to each inverter can simply be adjusted to take into account the effect of \jdw{$Z$}.

\subsection{Closed-loop representation of the overall system} \label{cl}
In this section we discuss how the power system can be viewed as the closed-loop interconnection of dynamical systems representing the bus dynamics and the \jdw{line} dynamics respectively. The bus dynamics are the dynamics at each bus $j \in N$ which relate the aggregate injected current $i^j_{DQ}$ to the bus voltage $v^j_{DQ}$. If an appropriate inverter is present, these dynamics will include the inverter dynamics given in Section~\ref{gfc}.

\begin{figure}[!t]
    \centering
    \includegraphics[width=0.3\textwidth]{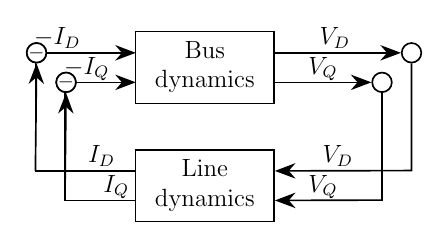}
    \caption{The power system as an interconnection of bus and line dynamics.}
    \label{fig:iosys}
\end{figure}

More precisely, the \jdw{line} dynamical system has input $V_{DQ}$ and output $I_{DQ}$ and dynamics given by~\eqref{netwdyn}. Note (as stated in Section~\ref{nw}) that $I_{DQ}$ is a vector of the total $DQ$ currents flowing into the network from each bus and $V_{DQ}$ is a similarly ordered vector of the $DQ$ bus voltages. We close the loop by considering bus dynamics with input $-I_{DQ}$ (the negative sign comes from the fact that this is current leaving the network and flowing into the bus) and output $V_{DQ}$. Fig.~\ref{fig:iosys} shows the representation of the power system as a negative feedback interconnection of the bus and line dynamics~\cite{spanias2018}.
For a bus $j \in N$ at which an inverter is present, Fig.~\ref{fig:iosys} and Fig.~\ref{fig:iofilt} are directly related as follows. The input $-i^j_{DQ}$ (which is a component of the vector $-I_{DQ}$ as in Fig.~\ref{fig:iosys}) is the current $-i_{DQ}$ in Fig.~\ref{fig:iofilt}. Similarly, the voltage $v^j_{DQ}$ (also a component of $V_{DQ}$ in Fig.~\ref{fig:iosys}) is the inverter capacitor voltage $v_{DQ}$ in Fig.~\ref{fig:iofilt}. Therefore every bus has a two-input two-output system associated with it and the direct sum of these systems is denoted as bus dynamics in Fig.~\ref{fig:iosys}, i.e. it relates $V_{DQ}$ with $I_{DQ}$ and has a corresponding block diagonal structure.

\subsection{Passivity}\label{passsect}
Passivity approaches have been widely used in the literature to achieve network stability in a decentralized and scalable way. Since the negative feedback interconnection of two passive systems is stable and passive~\cite{khalil1991}, by representing the system dynamics as an interconnection of two subsystems and ensuring that both are passive with respect to their inputs and outputs, we can guarantee stability. We consider for each bus $j \in N$ a dynamical system of the form:
\begin{equation}\label{eq:ss}
\begin{matrix}
    \dot{x_j} = f_j(x_j, w_j)\\
    z_j = g_j(x_j, w_j)
\end{matrix}
\end{equation}
where $f_j$ \il{is locally Lipschitz and $g_j$ is continuous}. \jdwrev{The inverter model and line model in this paper are linear, however bus dynamics of some loads or various inverter controllers may be nonlinear and hence we present the bus dynamics in the general form of~\eqref{eq:ss}.}

\begin{defn}[Passivity, \jdwrev{\cite[Def 6.3]{khalil1991}. }] \label{passdefn}
\jd{A system as in~\eqref{eq:ss} is said to be locally passive about some equilibrium point $(w_j^*, x_j^*)$ if there exist open neighbourhoods $W_j$ of $w_j^*$ and $X_j$ of $x_j^*$ and a continuously differentiable, positive semi-definite function $\jw{\mathcal{V}}_j(x_j)$ such that, for all $w_j \in W_j$ and all $x_j \in X_j$:}
\begin{equation}\label{eq:assu3}
    \dot{\mathcal{V}}_j \leq (w_j-w_j^*)^T(z_j - z_j^*)
\end{equation}
\jdw{\il{If the property holds with the right hand side of~\eqref{eq:assu3} replaced by $(w_j-w_j^*)^T(z_j - z_j^*) - (x_j - x_j^*)^T\phi_j(x_j - x_j^*)$, where $\phi_j$ is a positive definite function, then the system is said to be locally strictly passive about $(w_j^*, x_j^*)$. If instead the right hand side of~\eqref{eq:assu3} \ic{is} replaced by $(w_j-w_j^*)^T(z_j - z_j^*) - (z_j - z_j^*)^T\phi_j(z_j - z_j^*)$, then the system is said to be locally output-strictly passive about $(w_j^*, x_j^*)$.}} \jwr{We say that an output-strictly passive \li{system 
has output-strict passivity} index \kkl{$\rho_j>0$ }if \li{condition \eqref{eq:assu3} is replaced by}:
\begin{equation}\label{eq:osp}
   \kkl{ \dot{\mathcal{V}}_j \leq	(w_j-w_j^*)^T(z_j - z_j^*) - (z_j - z_j^*)^T\rho_j(z_j - z_j^*)}
\end{equation}}
\end{defn}
\begin{rmk}
\jdw{It should be noted that if the system~\eqref{eq:ss} is linear then if the passivity property holds locally about an equilibrium point, it also holds globally, i.e., for any other equilibrium point and for any deviation from an equilibrium point. Therefore for linear systems which satisfy Definition~\ref{passdefn}, we will refer to them as being "passive"
\il{(or "strictly passive")}
without making reference to a particular equilibrium point.}
\end{rmk}

\ilc{We start} by showing \ilc{in Proposition~\ref{passline}} that the line dynamical model is passive in \ic{a} common reference frame, \ic{irrespective of its complexity, i.e. for an arbitrarily large number of elementary components used to represent the line dynamics and despite the parameter variation along the lines}.
\jdw{This \ic{follows from the fact that the lines are composed of passive RLC components} \il{and the DQ transform on a common reference frame is an orthogonal transform and thus preserves \ilc{passivity}}.} 
\jdwrev{The associated proof may be found in Appendix B}.
\begin{proposition}
[Passivity of line dynamical model]\label{passline}
The model of the network lines~\eqref{netwdyn} \il{with} input $V_{DQ}-V^*_{DQ}$ \il{and} output $I_{DQ}-I^*_{DQ}$ \il{is \il{strictly} passive}.
\end{proposition}

\begin{rmk}
\jdw{Since the line dynamics are linear, the passivity property holds for any equilibrium point $(V^*_{DQ}, I^{N\ast}_{DQ})$.} However, Proposition~\ref{passline} does not generally hold when using local $dq$ coordinates due to the different frequencies and angles at various buses. This is a key reason for using a common reference frame for control of every inverter in the network.
\end{rmk}

Given the result in Proposition~\ref{passline}, passivity allows us to guarantee the stability of the general power network if all the bus dynamics are passive in the common reference frame. We state this result formally in \ilc{Theorem~\ref{stab} and the proof is provided in Appendix C.}

\begin{assumption}\label{ospass}
We assume that the bus dynamical systems satisfy an \jdw{output-strict} passivity property of Definition~\ref{passdefn} about an equilibrium $(w_j^*, x_j^*)$ for each bus $j \in N$. The storage function $\mathcal{V}_j(x_j)$ is continuous and has a strict local minimum at $x_j^\ast$. We further assume that for a constant input $w_j^*$, the equilibrium $x^*_j$ is locally \ilc{asymptotically stable.}
\end{assumption}
\jwr{Assumption~\ref{ospass} presents a passivity condition that must be satisfied at all buses, including both inverter buses and load buses. The following remarks discuss the significance of Assumption \ref{ospass} for both generation (inverter) and load buses.}
\begin{rmk} \label{genrmk}
\jwr{The passivity condition in Assumption~\ref{ospass} must be satisfied at all buses, and our focus is on the generation / inverter buses. We therefore aim to design grid-forming inverter control schemes which satisfy this Assumption, leading to the stability result we will shortly present. In addition, we will discuss how an appropriate control policy can additionally satisfy other conditions for power sharing or performance.}
\end{rmk}

\begin{rmk} \label{loadrmk}
\jwr{Assumption~\ref{ospass} also has implications for the incorporation of loads in our framework. If the load dynamics \li{satisfy this Assumption, i.e. they are passive about the equilibrium point, they can be incorporated in our analysis.} As it is well known, \jjwr{this property is satisfied by a range of loads including} constant impedance \jwr{loads} (Z) and some non-linear loads under appropriate conditions~\cite{strehle2020}. Individual non-passive loads may also be passivated via strictly passive load at the same bus such that the aggregate dynamics are passive. \jwr{A further possibility is to \kl{incorporate the load dynamics with the inverter dynamics at the bus \jwr{(if such an inverter is present), with} the objective of passivating the aggregate bus dynamics by an appropriate inverter control design. \jwr{The strictness of the passivity property achieved using the control schemes in this paper makes this \jjwr{possible}.}}}}
\end{rmk}


\begin{thm}[Stability]\label{stab}
\jdw{Suppose there exists an equilibrium $\jwrev{q^\ast} = (I^{N*}_{DQ}, x_1^*, ..., x_{|N|}^*)$ of the interconnected bus dynamics and line dynamics~\eqref{netwdyn}, for which the bus dynamics satisfy Assumption~\ref{ospass} for all buses $j \in N$. Then such an equilibrium is asymptotically stable.}
\end{thm}




\begin{rmk}The advantage of the stability criterion in Assumption~\ref{ospass} is \ilc{the fact} that it is a fully decentralized condition which allows \jdw{a wide range of heterogeneous} dynamics at each bus, so long as the aggregate dynamics satisfy the passivity condition. Instead of requiring the detailed modeling of the entire system, each bus may have \icl{different} dynamics which individually satisfy the passivity condition. This is therefore a decentralized and scalable approach.\end{rmk}

\section{Controller design} \label{design}

\jdw{\il{Even though a passivity condition in a common reference frame is a sufficient decentralized condition for stability, satisfying this property can in general be a non-trivial problem when higher-order bus dynamics are present. This is, for example, reflected in the fact that no simpler droop control scheme for grid-forming inverters in the literature satisfies a decentralized stability condition \iclc{that leads to a} plug and \jwwww{play} capability without requiring simplifications in the inverter or line dynamics.}} In this section we present our proposed control design
\iclc{such that} the passivity condition \iclc{is satisfied} along with various other objectives. We first explain the broader framework of the microgrid operation, list our control objectives, and then \jwwww{present} the control design.
\subsection{Microgrid operation \jjdw{and control objectives}} \label{design_mg}

\jwrev{A} centralized PMS is used to coordinate the inverters, as in other non-droop approaches~\cite{sadabadi2017}. The PMS solves the optimal power flow to calculate the appropriate setpoints, which are broadcast to each inverter in the microgrid. Frequency control is open loop and is achieved by means of internal oscillators at the \ilc{nominal frequency.} 
\jjwr{Since small errors in frequency cause slow drifts and eventual loss of synchronism, in practice an appropriate communication protocol (e.g. using GPS) is used to synchronize these frequencies at a slower timescale.} Each inverter therefore operates at the same frequency and is controlled in the same synchronous $DQ$ frame of reference. If no new setpoint is sent to the inverter (e.g. due to communication failure), it will operate autonomously using the embedded virtual impedance approach. This allows our proposed approach to be reliable and flexible. 

\kl{Given this \jwwww{microgrid} context, our control objectives are as follows.
\begin{itemize}
	\item Ensure some performance requirements. These specifications typically include: voltage reference tracking, \li{restricted} control inputs and \jjwr{limited} coupling between buses.
	 \item Allow power requirements of unexpected load changes to be shared between multiple inverters.
	\item Ensure stability for a general network with plug-and-play functionality. This can be \jwwj{guaranteed} via Theorem~\ref{stab} by ensuring the \jjwr{passivity of the bus dynamics}. 
\end{itemize}}


\kl{\begin{rmk}{
			The significance of {the}  passivity constraint is that \jjwr{adding} this
			{to optimization problems associated with control design \jjwr{leads to control policies which provide} stability guarantees for the network.}}
\end{rmk}}

\kl{\begin{rmk}\label{rmq:L2_pass}
		\jwjw{\jjwr{Output-strict passivity is a way to} \ic{have performance guarantees} in the presence of \ic{disturbances} (e.g. load fluctuations) 
			\il{by means of a restricted}
			$L_2$-gain for each bus, e.g.~\cite[Theorem 2.2.15]{vanderschaft2016}. 
			\il{This is a property that can be ensured with output-strict passivity, and is hence
				one of the reasons for} \ic{maximizing output-strictness (parameter $\rho$) in the inverter control design is a reasonable performance criterion. 
				Additional performance criteria will be discussed \jjwr{in \li{Sections~\jwwj{\ref{sec:mixed_controller}, \ref{casestudy}}} and validated with simulations}}.}
\end{rmk}}

\subsection{Full-state feedback design}\label{fsfdesign}

The full-state feedback  for the inverter model~\eqref{eq:convdyn2} is of the form
$ u = -K\tilde{x} - Mw$
where $K\in\mathbf{R}^{2\times 6}$ and $M\in\mathbf{R}^{2\times 2}$ are gain matrices: two $dq$ control inputs \il{each being a function of the }
four states of the physical system, two states corresponding to the voltage integrators, and two input feedback terms. This state feedback requires sixteen control gains to be chosen.

Given the passivity property~\eqref{eq:osp},  Theorem~\ref{thmlmi} states a \jwrev{matrix inequality} through which this property can be satisfied (e.g.~\cite{mccourt2013}), \jdwrev{and the proof is given in Appendix D}. \jdwr{As we will show, these inequalities can be added as a constraint to the optimization problems associated with control design to ensure the passivity property (Assumption \ref{ospass}).} 

\begin{thm}\label{thmlmi}
\il{Consider the} state feedback controller $u = -K\tilde{x} - Mw$ \jwjw{in conjunction with the dynamics~\eqref{eq:convdyn2}.} \il{The system with input $w$ and output $z$ satisfies the} 
passivity condition in \li{Assumption~\ref{ospass}}
\jdw{\il{and is output-strictly passive with strictness index}}
$\rho$ \jdwr{if} there exists a symmetric matrix $P = P^T > 0$ that satisfies:
	\begin{equation}\label{lmi}	
    \begin{bmatrix}	
    A_c^TP+PA_c\kl{+}2\rho C_c^TC_c & PB_c - C_c^T \kl{+} 2\rho C_c^TD_c \\	
    B_c^TP - C_c \kl{+} 2\rho D_c^TC_c & 2\rho D_c^TD_c-D_c^T-D_c	
    \end{bmatrix}	
    \leq 0\\[2ex]
\end{equation}	
where $A_c = \tilde{A}-\tilde{B}_uK$, $B_c = \tilde{B}_w - \tilde{B}_uM$, $C_c = \tilde{C} - \tilde{D}_uK$, and $D_c = \tilde{D}_w - \tilde{D}_uM$. Note that in this case $\tilde{D}_u$ and $\tilde{D}_w$ are zero matrices.
\end{thm}
~\\
\begin{rmk}
\jdwr{As stated,~\eqref{lmi} is a bilinear matrix inequality. However, using the change of variables $Q = P^{-1}, Y = KP^{-1}$ and the Schur complement, it can be rearranged as an LMI. \jwr{We leave this in the \kl{BMI} form as it allows us to directly place bound on the control gains in $K$ and $M$, which is not straightforward if such a change of variables is applied. We will also discuss convex problems for \li{controller} synthesis in the next section. }}
\end{rmk}



\subsection{Mixed $H_\infty$/passivity controller  design}\label{sec:mixed_controller}
\kl{In order to achieve complex  closed-loop specifications for each inverter, we focus on the design of a dynamical controller~$\mathcal{K}$
	\begin{equation}
		\mathcal{K}: ~~~
		\begin{bmatrix}
			\dot{x}_K\\u
		\end{bmatrix}=
		\begin{bmatrix}
			A_K& B_K\\C_K& D_K
		\end{bmatrix}
		\begin{bmatrix}
			x_K\\y_{\text{control}}
		\end{bmatrix}
		\label{eq:dynamic_controller}
	\end{equation}	
	where  $x_K$ is the controller state vector, $A_K$, $B_K$, $C_K$ and $D_K$ are the controller matrices \li{in its state space representation},\ {$u$ is the control input as in \eqref{eq:convdyn2}}, and
	$y_{\text{control}}=[ 	(v^{ref}_{{DQ}})^T ~ -{i}_{DQ}^T ~ i_{iDQ}^T ~ v_{DQ}^T]^T$ represent the  vector of measured signals serving as inputs for the controller~$\mathcal{K}$. \jjwr{The virtual impedance is implemented via the reference $v^{ref}_{{DQ}} = v^{set}_{{DQ}} - Zi_{DQ}$.} This controller will be designed to ensure   the passivity from $-i_{DQ}$ to output $v_{DQ}$ and also other performance requirements   expressed as $H_\infty$ constraints to shape the frequency responses of \li{appropriate closed-loop transfer functions}:
	\begin{itemize}
		\item $T_{v^{ref}_{DQ} \rightarrow e_{DQ}}$ to obtain the desired voltage reference tracking requirements;
		\item $T_{v^{ref}_{DQ} \rightarrow u}$ to limit the control input magnitudes;
		\item $T_{-i_{DQ} \rightarrow v_{DQ}}$  and $T_{-i_{DQ} \rightarrow u}$ to ensure restricted coupling between buses.
\end{itemize} }

\kl{These control specifications are  taken into account   using an appropriate choice of weighting \li{filters
$W_e$,} $W_u$ and $W_d$ (\lic{will be discussed more extensively in} Section \ref{sec:HinfPass_case_study}) such that
	\begin{equation*}{\small \begin{matrix}
				&\left\|W_e  T_{v^{ref}_{DQ} \rightarrow e_{DQ}}\right\|_\infty\leq 1 & \left\|W_e  T_{-{i}_{DQ} \rightarrow v_{DQ}} W_d\right\|_\infty\leq 1
				\\[1ex]
				&\left\|W_uT_{v^{ref}_{DQ}  \rightarrow u}\right\|_\infty \leq 1~~~~ & \left\|W_u T_{-{i}_{DQ} \rightarrow u} W_d\right\|_\infty\leq 1~~~
		\end{matrix}  }
	\end{equation*}
	with   $\left\|T\right\|_\infty$ being the $H_\infty$ norm of the linear system $T$ and is defined as
	$${\small \left\|T\right\|_\infty= \sup_{\li{\omega\geq0}
}~~\overline{\sigma}\left( T(j\omega)\right) }.$$}

\kkl{Therefore,  we are interested in designing a   controller $\mathcal{K}$    ensuring  the output-strict passivity condition   from $-{i}_{DQ}$ to $v_{DQ}$ together with the  $H_\infty$  conditions.
	These advanced specifications can be formulated as a convex optimization problem with LMI constraints which can be solved efficiently using standard \li{solvers \cite{scherer1997}. This LMI formulation is stated in Appendix E (Theorem \ref{thmHinflmi}), and discussed further in Section~\ref{sec:HinfPass_case_study}}.}

\subsection{Controller implementation}
\begin{figure}[!ht]
    \centering
    \includegraphics[width=0.425\textwidth]{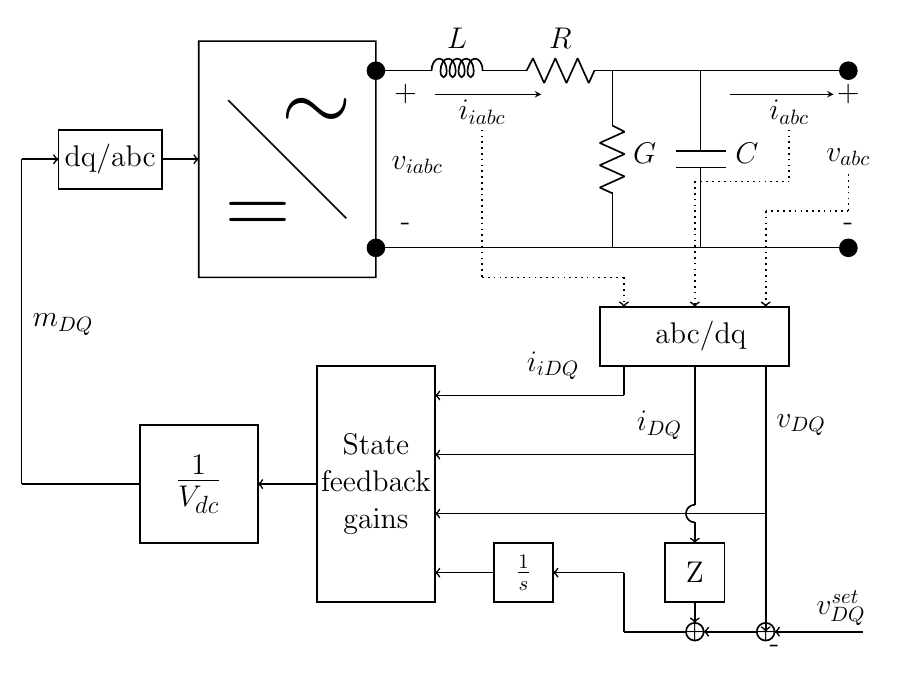}
    \caption{\jwrev{Block diagram of the 
    state feedback controller implementation and \li{system.}}}
    \label{fig:ctrlblk}
\end{figure}

\begin{figure}[!ht]
	\centering
	\includegraphics[width=0.9\linewidth]{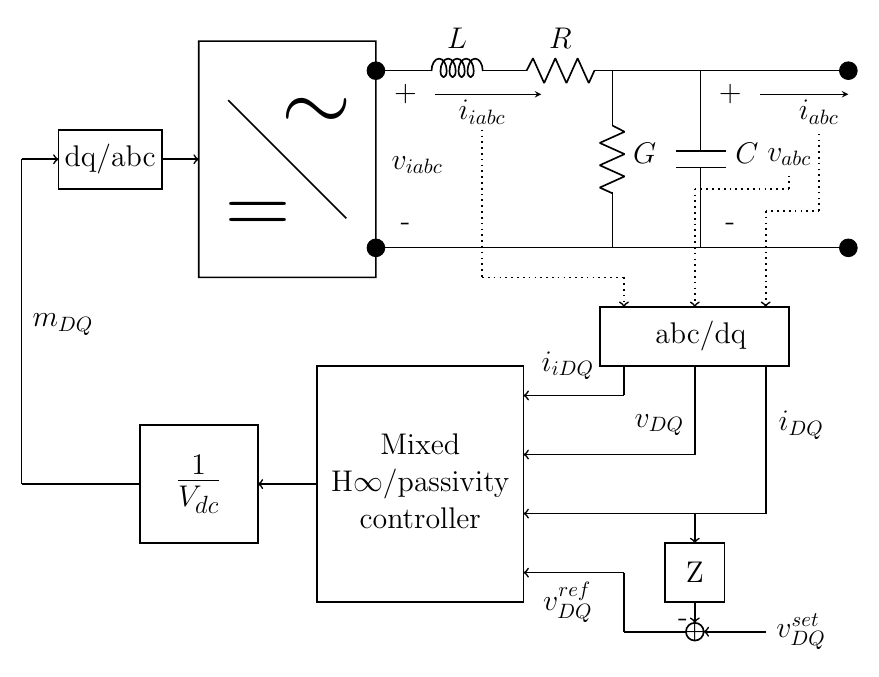}
    \caption{\jwrev{Block diagram \li{of the 
    $H_\infty$/passivity controller implementation and system. It should be noted that this is a dynamic controller.}}}
	\label{fig:dynctrlblockdiagram}
\end{figure}

\kl{We illustrate the \jdwr{full-state feedback} control design  and the mixed $H_\infty$/passivity control  design in Fig.~\ref{fig:ctrlblk} and  Fig.~\ref{fig:dynctrlblockdiagram}.} The full-state feedback design is effectively a single-loop design, in contrast to traditional double (inner current / outer voltage) loop design~\cite{chandorkar1993}. Single loop controllers are a promising alternative~\cite{wang2017} which have been successfully demonstrated, e.g.~\cite{qoria2019}. \jdw{The advantages of the single loop are simpler control design and better stability properties}~\cite{wang2017} which, in our case, can be exploited to guarantee stable plug-and-play operation. \jdwr{Furthermore, although full-state feedback designs are sometimes impractical due to the requirement to obtain full state information, in this case the only measurements required are the filter voltage and currents, which are also necessary for a traditional double PI loop controller.}

\jjwr{The major differences between the controllers in Figs.~\ref{fig:ctrlblk} and~\ref{fig:dynctrlblockdiagram} are as follows.  The state feedback controller is simple and \li{this allows numerically to consider in control synthesis also certain practical constraints which are non-convex.} By contrast, the mixed $H_\infty$/passivity controller \jjwr{is dynamic and of high order. \li{A main feature of this controller is that it allows 
$H_\infty$ criteria to be satisfied by formulating a} convex LMI problem.}}
\jdwrev{As always, there are \li{underlying trade-offs:} between \li{the simplicity of a static state feedback controller} and the flexibility of a high-order \li{dynamic} controller; and also between convexity of the optimization problem and the incorporation of certain non-convex constraints which may be useful in practice.}
\begin{rmk}
An issue when eschewing classical cascaded PI loops is current limitation \jdwrev{for protection}, since in most practical designs the current is limited via saturation of the current reference. 
\iclc{Virtual} impedance is a proven alternative~\cite{paquette2015} to limit the current and \jwwww{this can be achieved with our design}.
\end{rmk}
\begin{rmk}
\jw{Angle droop and virtual impedance are closely related~\cite{sun2017}. Compared to angle droop, our approach and control architecture offers several advantages. \jwwww{Firstly,} passivity \jwwww{(and thus plug-and-play operation)} is guaranteed, whereas angle droop schemes \jwwww{are only passive for small droop gains} \ilc{\cite{watson2019}.}
\jwwww{Secondly,} droop \iclc{approaches} limit the achievable performance due to the presence of only a few tunable parameters.}
\end{rmk}

\section{Case studies}  \label{casestudy}

\subsection{Control specifications} \label{cs-ctrl}

Although stability is guaranteed via the passivity condition, \jwjw{(likewise $L_2$ performance by the output-strictness), in practice further specifications are required to achieve good performance. A sensible approach, and the one we adopt, is to ensure good performance in the case of one inverter with a load, and then add constraints to \jwrev{limit output voltage deviations which also has the effect to reduce} the coupling between inverters. Although this method does not give general performance guarantees for any network, our investigations indicate that in practice appropriate network-wide transient performance is achieved.} 

\kl{In addition to the passivity constraint, the different performance specifications \jdwr{considered} in this paper are as follows.
\begin{itemize}
	\item Fast convergence of the voltage tracking error.
	\item Small control input magnitudes.
	\item \jjwr{Limited} coupling between inverters \jdwr{(in the form of loop shaping design in the frequency domain)}.
\end{itemize}}

\kl{These control specifications will be incorporated as constraints in the optimization problem formulation to obtain the desired controllers.}
\subsection{State feedback controller synthesis}
We start with the static state and input feedback controller $u=-K\tilde{x}-M w$ and look to design appropriate $K$ and $M$. \jdwr{The states of the plant $\tilde{x}$ include four states associated with the inverter filter and two states associated with the integrator.} The previous control specifications can be expressed as follows.
\jd{\begin{itemize}
    \item Maximum controller gains: Note that in general some limit of this sort is necessary to prevent a solution with impractically large \iclc{gains.}
    \item Maximum real part of eigenvalues \iclc{of $A_c$}: This determines the speed of the response for each inverter. \jww{The speed of the response is chosen to be sufficiently fast, \iclc{but} 
        must be limited to avoid a controller that saturates the inverter's PWM operation.}
    \item \iclc{Bound on the frequency response}
    \iclc{of the bus dynamics with input $-i_{DQ}$ and output $v_{DQ}$:}
    \iclc{This reduces the coupling between inverters.}
    %
\end{itemize}}
\iclc{It should also be noted that the virtual impedance is also a parameter \jdwrev{which can be chosen} (discussed previously) and provides a trade-off between voltage regulation and power sharing.}
\jdw{The optimization problem \jwwww{may be stated} as follows, where \jjjw{$\tilde{n}$} is the number of states\footnote{\jdwr{In this example $\tilde{n}$ is 6. In the case of a different inverter filter, or the inclusion of other DC or AC dynamics, $\tilde{n}$ may vary. However, since this is the number of states at only one bus it should generally remain small.}}}
\jdw{\begin{equation}\label{opt}
\begin{array}{c}
    \jww{\underset{K,M,P}{\max}}\phantom{l} \rho,  \text{ such that}\\
        \ilcl{\eqref{lmi} \text{\il{\mbox{ }holds}}}\\
    \jdwr{|K_{ik}|} \leq p_{max} \\
    \jdwr{|M_{ik}|} \leq p_{max} \\
    \Re[\lambda_i(A_c)] \leq \lambda_{max}\\
    \bar{\sigma}(\tilde{C}(j\omega I-A_C)^{-1}B_c+\tilde{D}_w) \leq |\gamma\frac{\omega_c}{j\omega+\omega_c}| \text{ for all }\omega
\end{array}
\end{equation}}
~\\

 \begin{table}[ht]
    \centering
    \caption{Tuning specifications \kl{for the state feedback controller}.}
    \begin{tabular}{c|c}
     Parameter & Value \\
     \hline
     Maximum controller \iclc{gains}
      & $p_{max}=125$ \\
     Maximum real part of eigenvalue & $\lambda_{max}=-5$ \\
     Frequency response (max gain, bandwidth) & $\gamma = 1.5$,  $\omega_c = 10^5$ rad/s\\
     Virtual impedance &  $R_V$ = 0.5, $X_V$ = 1 
    \end{tabular}
    \label{paramsCtrl}
\end{table}
\vspace{-3mm}

We solve the optimization problem numerically, using MATLAB's $systune$~\cite{apkarian2014}. This is a non-convex problem, nonetheless the controller obtained satisfies the passivity condition as well as the design specifications discussed before, resulting in a good design.

Using the inverter filter parameters as in Table~\ref{tab:analyticalparams} \jwwj{and specifications as in Table~\ref{paramsCtrl}} \jdw{\jdwrev{\jwwj{results} in a maximized value of $\rho = 0.4000$}}, \li{and} the following controller is obtained:
\begin{equation} \label{gains}
\begin{matrix}
    K =  \begin{bmatrix}117.3 & 1.1 & 6.3 & 0.4 & 40.0 & -7.3\\
    -2.6 & 117.2 & -2.1 & 12.9 & 2.1 & 72.5\end{bmatrix},\\
    \iclcl{M =}  \begin{bmatrix}107.8 & 3.3\\-1.2 & 104.7\end{bmatrix}
\end{matrix}
\end{equation}
\jwrev{In our experience, the problem is generally feasible for reasonable parameter choices. However, there are obvious cases where the problem may not be feasible; for example, specifying a very fast response while bounding the controller gains $p_{max}$ to a low number; or bounding the maximum gain of the frequency response too low. It is also not a feasible problem to have a negative virtual resistance as well as passivity since such constraints are not compatible. Should the problem be infeasible, it is recommended to adjust the parameters (e.g. higher $p_{max}$) or set less ambitious specifications.}

\kl{\begin{rmk}
		The  state feedback problem of~\eqref{opt} is non-convex as it involves \li{a BMI and} non-convex constraints.  Nevertheless, \jdwr{it} should also be noted that the  state feedback problem is local and thus always small in dimension,  hence the non-convex problem is computationally tractable, i.e. controllers with good performance can be obtained in a fraction of a second on a modern  computer. 
\end{rmk} }

\subsection{\kl{Mixed $H_\infty$/passivity  controller synthesis}}\label{sec:HinfPass_case_study}
\kkl{We consider now the controller~\eqref{eq:dynamic_controller}. The objective is to compute the \lic{matrices}~$A_K$, $B_K$, $C_K$ and $D_K$ such that  the $H_\infty $ performance constraint together with passivity are satisfied.}

\kkl{First, {we} consider the performance requirements. As mentioned earlier and as {usually} done \li{in $H_\infty$ control,} the performance requirements are specified via the weighting filters   $W_e$, $W_u$ and $W_d$   \li{as follows:}
	\begin{itemize}
		\item    Fast   convergence of the  voltage tracking error $e_{DQ}=v^{ref}_{DQ}-v_{DQ}$ with a minimal steady-state error. This can be \li{achieved} by shaping \li{$T_{v^{ref}_{DQ} \rightarrow e_{DQ}}(j\omega)$}
using $W_e(j\omega)$ such that \li{for all $\omega\geq0$}			
		$${ \overline{\sigma} \left( T_{v^{ref}_{DQ} \rightarrow e_{DQ}}(j\omega)\right)\leq\overline{\sigma} \left( W_e(j\omega)^{-1}\right)  }$$
		with
		\begin{equation}
			W_e(s)=\text{diag}\left(\dfrac{M^{-1}_{s_1}s+\omega_{b_1}}{s+\epsilon \omega_{b_1}}, \dfrac{M^{-1}_{s_2}s+\omega_{b_2}}{s+\epsilon \omega_{b_2}}\right)
			\label{eq:W_e}
		\end{equation}
		where the parameters $\omega_{b_i}$, $M_{s_i}$ and $\epsilon$ allow to impose the desired closed-loop bandwidth, the minimum peak value and the maximum allowed steady-state error respectively.
		\item Allowed control input magnitude by imposing an upper bound on \li{$\overline{\sigma}\left( T_{v^{ref}_{DQ} \rightarrow u}(j\omega)\right)$} such as
		$$\overline{\sigma} \left( T_{v^{ref}_{DQ} \rightarrow u}(j\omega)\right)\leq
\li{ W_u^{-1}}
$$
		\li{for all $\omega\geq0$,} with 	\begin{equation}
\li{	W_u=\dfrac{1}{M_{u_1}}}
\label{eq:W_u}
		\end{equation} \jjwr{where \li{$M_{u_1}$} 
is a \li{constant} 
that limits the control input gain.}
		\item Rejection of high frequency components of $-i_{DQ}$ and their effect on  $v_{DQ}$ and $u$.  This can be done by shaping $  T_{-{i}_{DQ} \rightarrow v_{DQ}}$  and $ T_{-{i}_{DQ} \rightarrow u}$ using $W_d(j\omega)$ such that	\li{for all $\omega\geq0$}
		$$\overline{\sigma} \left( T_{-{i}_{DQ} \rightarrow v_{DQ}} (j\omega)\right)\leq\overline{\sigma} \left( W_d(j\omega)^{-1}\right)  $$
		$$\overline{\sigma} \left( T_{-{i}_{DQ} \rightarrow u}(j\omega)\right)\leq\overline{\sigma} \left( W_d(j\omega)^{-1}\right)  $$			
		with
		\begin{equation}
			W_d(s)=\text{diag}\left(\dfrac{s+\omega_{d_1}}{\varepsilon s +\omega_{d_1}},\dfrac{s+\omega_{d_2}}{\varepsilon s+ \omega_{d_2}} \right) 				\label{eq:W_d}
		\end{equation}
		where the parameters $\omega_{d_i}$ \li{and $\varepsilon$} allow to impose the desired cut-off frequency   and the maximum allowed magnitude at high frequencies.		
\end{itemize} 	}

\kkl{All \jwr{these} frequency loop shaping constraints are satisfied~if  \begin{equation}\left\|\begin{matrix}
			W_e  T_{v^{ref}_{DQ} \rightarrow e_{DQ}} & W_e  T_{-{i}_{DQ} \rightarrow v_{DQ}} W_d
			\\W_uT_{v^{ref}_{DQ}  \rightarrow u}  & W_u T_{-{i}_{DQ} \rightarrow u} W_d
		\end{matrix} \right\|_\infty \leq 1
		\label{eq:Hinf_con}
\end{equation}}

\kkl{Now, let $\eta$ be the output-strict passivity index  of $T_{-{i}_{DQ} \rightarrow v_{DQ}}$ which is the transfer function from  $-{i}_{DQ}$ to $v_{DQ}$ (including weighting filters).
	Similarly to \jjwr{Section~IV.B}, \li{the 
following optimization problem can be considered} 
\begin{equation*}
		\begin{array}{l}
			\max_{}    \eta \\
			~~~\text{s.t. \li{the following hold}}\\
			~~~~~\text{Output-strict passivity of $T_{-{i}_{DQ} \rightarrow v_{DQ}}$ as in~\eqref{eq:osp}}\\	~~~~~\text{$H_\infty$ condition in~\eqref{eq:Hinf_con}}\\
		\end{array}	
\end{equation*}}

\kkl{
	Using the results in \cite{scherer1997},    a convex optimization problem  with  LMI constraints can be formulated  that leads to a controller satisfying both the output-strict passivity of $T_{-{i}_{DQ}}  \rightarrow v_{DQ}$ as in~\eqref{eq:osp} and the $H_\infty$ condition in~\eqref{eq:Hinf_con}. The LMI formulation is \jjwr{stated} in Theorem 3 in \jjwr{Appendix} E.}

\kkl{Thereafter, using  the inverter filter parameters as in \jjwr{Table~\ref{paramsCtrl}}  and specifications as in \jjwr{Table~\ref{tab:paramsCtrl_Hinf}}, \li{maximizing} $\eta$ (or equivalently \li{minimizing} $\eta^{-1}$) \li{using the LMI formulation in Theorem~3} in Appendix E    results in a maximized value of~$\eta = 0.0132$.}
\kl{\begin{table}[ht]
		\centering
		\caption{\kl{Tuning specifications for the mixed $H_\infty$/passivity controller}.}
		\begin{tabular}{l|l}
			Parameter & Value \\
			\hline
			Desired closed-loop bandwidth of $T_{v^{ref}_{D} \rightarrow e_{D}}$	& $\omega_{b_1}=10$ rad/s\\
			Maximum peak value of $T_{v^{ref}_{Q} \rightarrow e_{Q}}(\mathbf{j}\omega)$ 	& $M_{s_2}=2$   \\	
			Desired closed-loop bandwidth of $T_{v^{ref}_{D} \rightarrow e_{D}}$	& $\omega_{b_1}=10$ rad/s\\
			Maximum peak value of $T_{v^{ref}_{Q} \rightarrow e_{Q}}(\mathbf{j}\omega)$ 	& $M_{s_2}=2$   \\			
			Maximum steady-state error & $\epsilon=10^{-3}$\\
			\hline
			Maximum control input gain & $M_{u_1}=300$   \\
			\hline
			Desired closed-loop cutoff frequency of $T_{i_{D} \rightarrow v_{D}}$	& $\omega_{d_1}=10^{4}$ rad/s\\	
			Desired closed-loop cutoff frequency of $T_{i_{Q} \rightarrow v_{Q}}$	& $\omega_{d_2}=10^{4}$ rad/s\\	
			\begin{tabular}{l}
				\hspace{-0.3cm}	Maximum allowed high frequency magnitudes\\ of $T_{i_{DQ} \rightarrow v_{DQ}}$
			\end{tabular}
			& $\varepsilon=10^{-3}$   \\				
		\end{tabular}
		\label{tab:paramsCtrl_Hinf}
	\end{table}
}

\subsection{Numerical \ilc{simulations}} \label{cs-num}

\begin{figure}[!ht]
    \centering
    \vspace{-3mm}
    \includegraphics[width=0.4\textwidth]{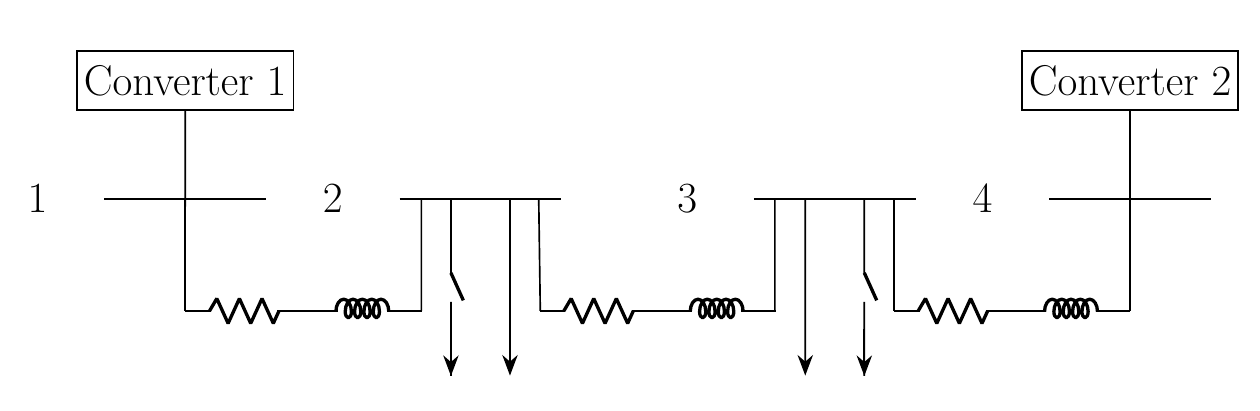}
    \vspace{-3mm}
    \caption{Diagram of microgrid used for simulation study.}
    \label{fig:testnet}
\end{figure}

In this section we verify the passivity approach to grid-forming converter design via \icl{simulations conducted using Simscape Power Systems}. The simulation model is detailed and realistic, and includes the PWM switching of the inverters and the DC-side dynamics, switching transients, etc. We consider a four-bus test network with two converter and load buses as shown in Fig.~\ref{fig:testnet}. The network parameters \jjww{are given} in Table~\ref{tab:params}. \jjww{For convenience \ilc{identical inverters are used in the simulation model.} \kl{The proposed  {state feedback} and $H_\infty$/passivity \jwr{controllers} have the design specifications given in Table~\ref{paramsCtrl}  and Table~\ref{tab:paramsCtrl_Hinf}  respectively;} and \jdwr{are} compared to other control schemes as in {Table~\ref{tab:analyticalparams}}, where the inverter filter parameters are also presented.}

\jwrev{The static loads are modeled as constant impedance, while the switched loads are modelled as part (half) constant power, part constant impedance}.
Our case study is comprised of \jwww{two} parts: firstly, a comparison of the transient response to a load change \il{with} 
existing alternatives, namely frequency, angle droop, and matching control schemes; secondly, \jjww{a demonstration of plug-and-play operation}.

\begin{table}[ht]
    \centering
    \caption{Network parameters.}
    \begin{tabular}{c|c}
     Description & Value \\
     \hline
     Static loads & 3 kW / 0.5 kVAr at buses 2 \& 3\\
     Switched loads & 4.5 kW / 0.5 kVAr at buses 2 \& 3\\
     Line parameters (1,2), (3,4) & 0.1 $\Omega$, 0.6 mH \\
     Line parameters (2,3) & 0.1 $\Omega$, 5 mH \\
    \end{tabular}
    \vspace{-3mm}
    \label{tab:params}
\end{table}

\begin{table}[ht]
    \centering
    \caption{Inverter parameters.}
    \begin{tabular}{c|c|c}
     Category & Description & Value \\
     \hline
     & Inverter filter ($R$, $L$, $G$, $C$) & (0.1 $\Omega$, 8 mH, 1/350 S, 50 $\mu$F) \\
     \hline
     & Frequency droop gain & $m_p$ = 9.4 $\times 10^{-5}$ \\
     & Angle droop gain & $m_p$ = 2 $\times 10^{-6}$ \\
     Droop & Voltage droop gain & $n_q$ = 1 $\times 10^{-4}$ \\
     control & Inverter voltage setpoint & $v_{DQ}^{set}$ = (311 V, 0 V) \\
     & Voltage loop PI gains & (5, 5) \\
     & Current loop PI gains & (5, 10) \\
     \hline
     & DC voltage, capacitance & 1000 V, 8 mF \\
     Matching & DC voltage loop PI gains & (1, 10) \\
     control & Other parameters & $\eta$ = 0.3142, $\mu$ = 0.311 \\
     & Parasitic conductances & $G_{dc} = G_{ac} = 0$ \\
     & Voltage droop gain & $n_q$ = 1 $\times 10^{-4} \times \frac{2}{V_{DC}}$
    \end{tabular}
    \label{tab:analyticalparams}
\end{table}

\begin{figure}[!ht]
    \centering
    \vspace{-3.5mm}
    \includegraphics[width=0.44\textwidth]{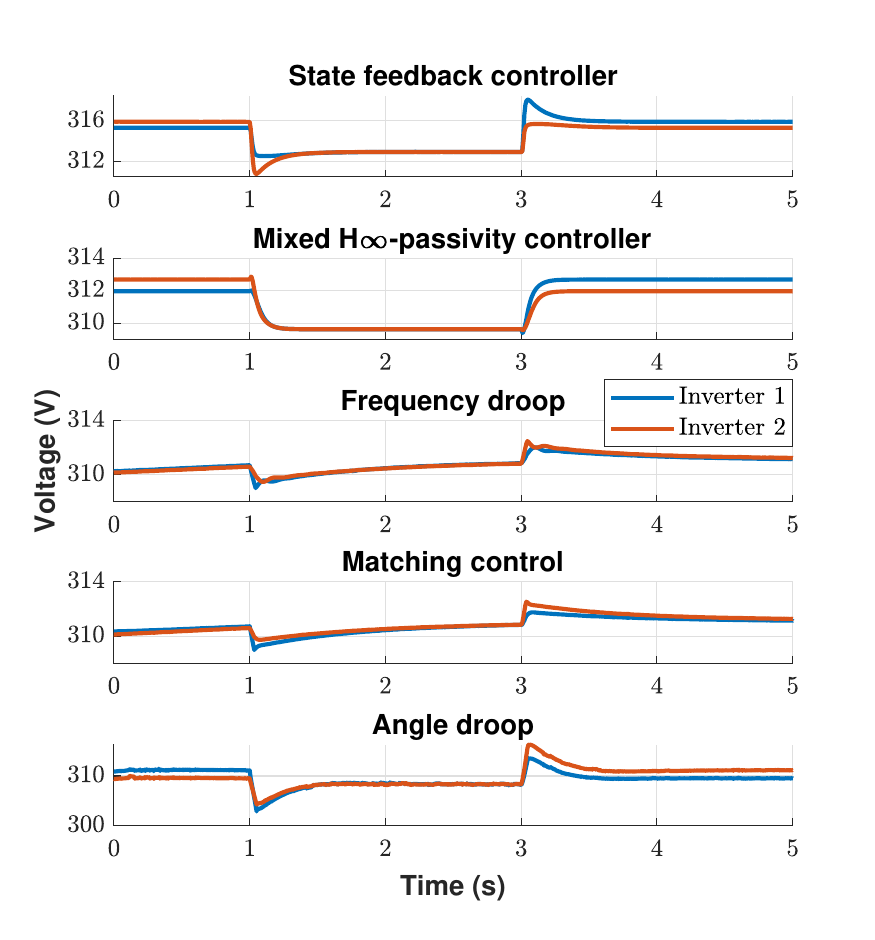}
    \vspace{-3mm}
    \caption{System response to load change: voltage.}
    \label{fig:voltresponses}
\end{figure}

\begin{figure}[!ht]
    \centering
    \vspace{-3mm}
    \includegraphics[width=0.429\textwidth]{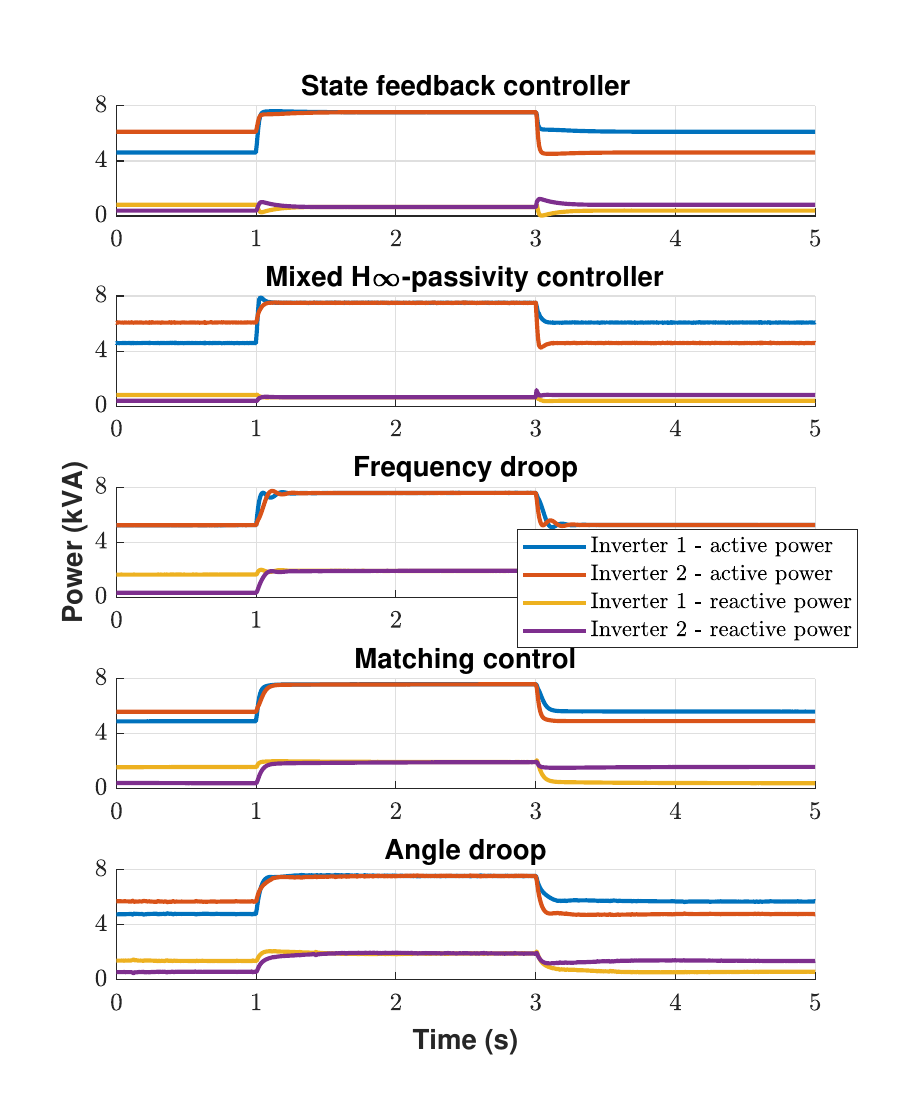}
    \vspace{-3mm}
    \caption{System response to load change: power.}
    \vspace{-3mm}
    \label{fig:powerresponses}
\end{figure}

In Figs.~\ref{fig:voltresponses}-\ref{fig:powerresponses} we show 
\iclcl{power} and voltage magnitude responses 
\ilc{to two step disturbances}: a 4.5 kW load added at bus 2 at t = 1s and the switching off of an equivalent load at bus 3 at t = 3s. In order to show that our approach can allow a 
\iclcl{faster response}
than existing alternatives, we also compare the response to \ilc{alternative control schemes} in the literature: frequency droop, angle droop,~\cite{majumder2009} and matching control schemes using parameter values from~\cite{chandorkar1993} for droop control and~\cite{arghir2018} for matching control as in Table~\ref{tab:analyticalparams}. \jd{A wide range of parameters were tried for the droop/matching schemes in order to achieve good behavior in the simulations~\cite{watson2019}.} \jdw{The appropriate equations of the three controllers (frequency droop, angle droop, and matching control) to which we compare our proposed \jdwr{schemes} may be found in our previous work~\cite{watson2019}, along with further discussion.} We note that the responses 
\iclc{of the proposed \jdwr{schemes}} are faster 
and the voltage regulation is \jjdw{appropriate. In practice, new voltage setpoints will be calculated and broadcast to correct any voltage deviations or power sharing inaccuracies.} 


\begin{figure}[!ht]
    \centering
    \includegraphics[width=0.429\textwidth]{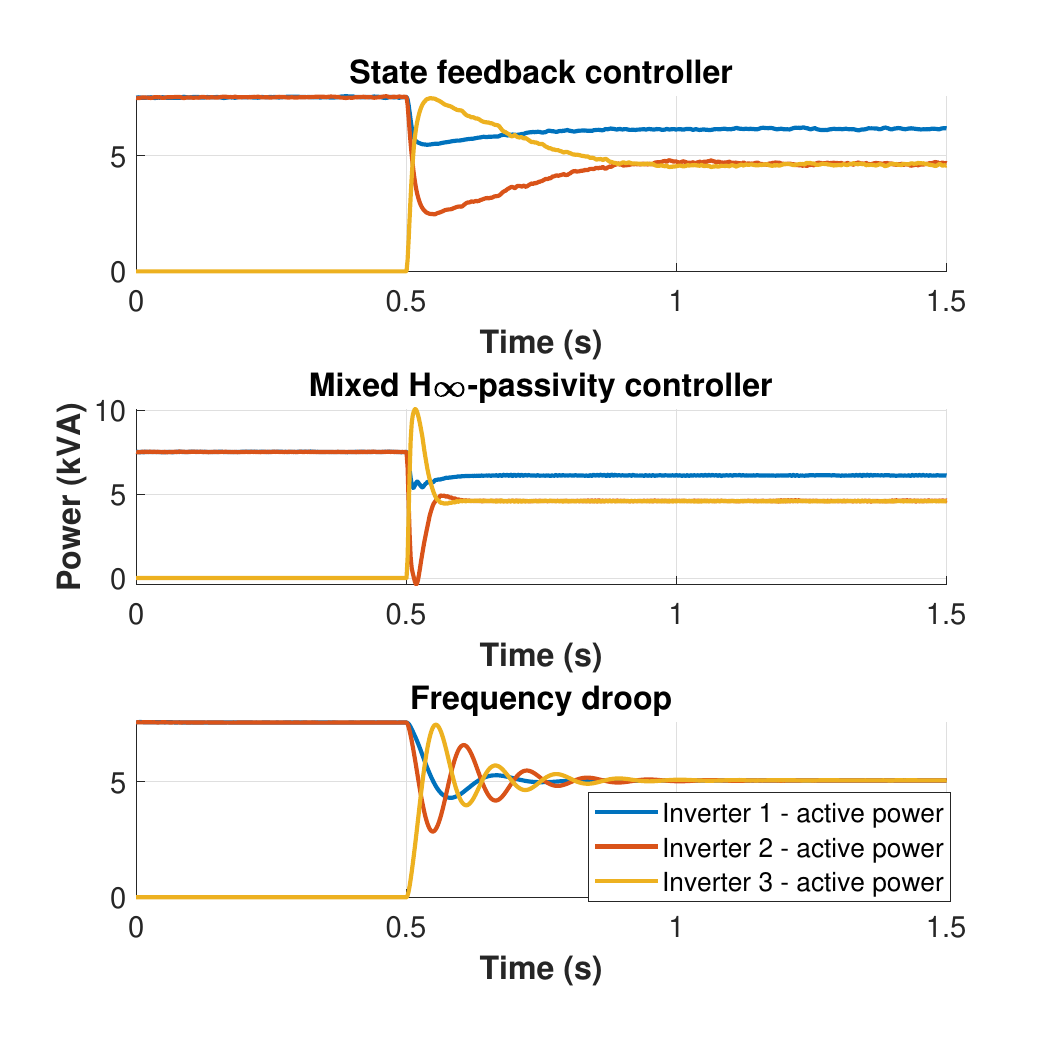}
    \vspace{-3mm}
    \caption{Plug-and-play response.}
    \label{fig:pppq}
    \vspace{-3mm}
\end{figure}

\jw{As discussed, one of the key advantages of the proposed \jdwr{control design framework} is that plug-and-play operation is possible. We demonstrate this \ilc{by simulating} the plug-in of an additional inverter at bus 4 in the test network (Fig.~\ref{fig:testnet}). The response \jdwr{using the state feedback \li{and the mixed $H_\infty$/passivity controllers}} is shown in Fig.~\ref{fig:pppq}, where this is also compared to traditional droop control. In the case of frequency droop, the phase angles of inverter 2 and 3 are synchronized exactly before connection, \jwjw{whereas in the proposed design no adjustment is necessary, and the inverters are unaware of the plug-in of another inverter. This is a challenging scenario that can lead to bad performance.} As illustrated, plug-and-play operation is possible with the proposed controller and less severe oscillations are experienced than in the case of the traditional frequency droop controller. }
\vspace{-3mm}
\section{Conclusion} \label{conclusion}
We have considered the problem of designing controllers for grid-forming converters in microgrids that regulate the voltage and frequency satisfactorily, achieve a power sharing property, and allow \ic{stability to be ensured via decentralized conditions.}
To this end, we propose a passivity-based controller design method for converter-based microgrids. Our proposed method achieves good performance while guaranteeing stability 
\ic{in a decentralized way}
via its passivity property. Furthermore, our controller improves the power sharing \icl{properties} 
of the microgrid compared to previous non-droop approaches. Finally, realistic simulations are used to verify \icl{the performance of} the proposed \jdwr{controllers}. A comparative study and a demonstration of plug-and-play operation 
were included to show that the controller design improves the stability and dynamic response of the microgrid in response to both load and setpoint changes.

\section*{Appendix}
In the appendices, we provide various definitions for the transmission line model as well as the proofs of all our results.
\subsection{Transmission line model}
\li{In this section we define the matrices $H$, $H_B$, $H_C$ that allow a concise representation of the line dynamics in~\eqref{eq:linemodel}. We start by defining various notions that are needed to specify these~matrices.}
\jwrev{
\begin{itemize}
    \item We consider a graph $(N_{L}, E_L)$ associated with the general RLC network representing the line dynamics in this case.
        In particular, we assign a node to each voltage in Fig.~\ref{fig:linemodel} (i.e. ($V^i_{dq}$, $V^{ij}_{dq,1}$, ... $V^{ij}_{dq,n-1}$, \jwrev{$V^j_{dq}$})) and do this for all lines. There exists an edge between any two consecutive nodes along a line.  Thus there is an edge associated with each current (($I^{ij}_{dq,1}$, $I^{ij}_{dq,2}, ..., I^{ij}_{dq,n-1}$) in Fig.~\ref{fig:linemodel}) and this is true for all lines.
        $N_L$ is therefore the set of nodes $(1,2, \dots |N_L|)$ representing the voltages of the buses and internal capacitances, and $E_L \subset N_L \times N_L$ is the set of edges associated with the current of each line subsection.
        A direction is assigned to each edge in $E_L$ which can be arbitrarily chosen, and similarly the ordering of the edge in $E_L$ can be arbitrarily chosen. The nodes in the set $N_L$ are ordered such that a node $j \in N_L, j \leq |N|$ corresponds to the bus $j \in N$.
    \item We enumerate the capacitances in the line model (for all lines), and the indices $(1, 2, ..., |N_C|)$ form the set $N_C$. The order of the enumeration corresponds to the order in $N_L$. To be precise, the node $i \in N_L, i = |N| + j$ corresponds to the internal capacitance $j \in N_C$.
    \item $\mathcal{R} > 0$, $\mathcal{L} > 0$, $\mathcal{G} > 0$, $\mathcal{C} > 0$ are diagonal matrices of the equivalent line resistances, inductances, conductances, and capacitances respectively with respective dimensions $|E_L| \times |E_L|$, $|E_L| \times |E_L|$, $|N_C| \times |N_C|$, $|N_C| \times |N_C|$. The order of $\mathcal{R}$ and $\mathcal{L}$ corresponds to the order of the edges in $E_L$, and the order of $\mathcal{G}$ and $\mathcal{C}$ corresponds to the order of the nodes representing internal capacitances in $N_C$.
    \item $I^\mathcal{N}_{dq}$ is the vector of length $2|E_L|$ which represents the current (i.e. $I^{ij}_{dq,k}$ for all $(i,j) \in E$ and corresponding $k$) of every edge in $E_L$. This vector is ordered consistently with the order of the edges in $E_L$.
    \item $V^C_{dq}$ is a vector of length $2|N_C|$ which represents the voltages of the capacitors (i.e. $V^{ij}_{dq,k}$ for all $(i,j) \in E$ and corresponding $k$ in the line representation~\eqref{eq:linemodel}), ordered correspondingly to the order in $N_C$.
\end{itemize}
We \ic{\lif{now}} define the matrices below with elements $0,1,-1$:
\begin{itemize}
    \item $H$ is the $|N_L|\times|E_L|$ incidence matrix of the graph ($N_L$, $E_L$) which is defined as:\\
    $
\begin{cases}
    H_{iz} = 1 &\text{if edge } z \text{ leaves node }i \in N_L\\
    H_{iz} = -1 &\text{if edge } z \text{ enters node }i \in N_L\\
    H_{iz} = 0 &\text{otherwise}
\end{cases}
$\\
i.e., $H$ illustrates for each edge in $E_L$ the nodes in $N_L$ it is connected to.
    \item $H_B$ is an $|N|\times|E_L|$ submatrix of $H$ that  consists of the first $|N|$ rows of $H$; i.e., it illustrates the connection between edges/nodes for edges connected to buses in $N$.
    \item $H_C$ is an $|N_C|\times|E_L|$ submatrix of $H$ that consists of the last $|N_C|$ rows of $H$; i.e., it illustrates the connection between edges/nodes for edges connected to nodes associated with capacitances in $N_C$.
\end{itemize}
Hence, $H_C$ and $H_B$ are related to the incidence matrix $H$ by means of the equation 
\kkl{$H = \begin{bmatrix}
H_B^T &
H_C^T
\end{bmatrix}^T$}.
}
\subsection{Proof of Proposition~\ref{passline}}
\jdwrev{This follows by using $\mathcal{V}_\mathcal{N}$ defined below to show that the condition for strict passivity stated in  Definition~\ref{passdefn} is satisfied.
\begin{equation}\label{eq:netwps}
\begin{aligned}
    \mathcal{V}_{\mathcal{N}} = \frac{1}{2}&(V_{DQ}-V^{\ast}_{DQ})^T [\mathcal{L} \otimes I_2](I^N_{DQ}-I^{N\ast}_{DQ}) \hspace{1cm}\\
    &+  \frac{1}{2}(V^C_{DQ}-V^{C\ast}_{DQ})^T [\mathcal{C} \otimes I_2](V^C_{DQ}-V^{C\ast}_{DQ})
    \end{aligned}
\end{equation}
By taking the time derivative and substituting~\eqref{netwdyn} we have
\begin{equation}\label{eq:storge_r}
\begin{aligned}
    \dot{V}_\mathcal{N} =& (V_{DQ}-V_{DQ}^*)^T(I^N_{DQ}-I^{N*}_{DQ})\\& - (x_\mathcal{N}-x_\mathcal{N}^*)^T{\small\begin{bmatrix} \mathcal{R} \otimes I_2 & 0 \\ 0 & \mathcal{G} \otimes I_2\end{bmatrix}}(x_\mathcal{N}-x_\mathcal{N}^*)
    \end{aligned}
\end{equation}
where $x_\mathcal{N} = [(I^{N}_{DQ})^T \mbox{ } (V^{C}_{DQ})^T]^T$. The matrix ${\small\begin{bmatrix} \mathcal{R} \otimes I_2 & 0 \\ 0 & \mathcal{G} \otimes I_2\end{bmatrix}}$ is positive definite, which completes the proof.} 

\subsection{Proof of Theorem~\ref{stab}}
\jdw{The result in Proposition \ref{passline}, together with equation \eqref{eq:assu3}, allows a straightforward LaSalle argument using the candidate Lyapunov function \il{$\mathcal{V}(q) = \mathcal{V}_\mathcal{N} + \sum_{j=1}^{|N|}\mathcal{V}_{j}(x_{j})$, where $q= ( \jwjwjw{x_\mathcal{N}}, x_1,x_2, ..., x_{|N|})$, \iclc{$x_\mathcal{N} = [(I^{N}_{DQ})^T \mbox{ } (V^{C}_{DQ})^T]^T$, and $\mathcal{V}_\mathcal{N}$ as in \eqref{eq:netwps}} to prove asymptotic stability. }}\jdw{Consider the time derivative of \il{$\mathcal{V}(q)$}. From \eqref{eq:storge_r}, we have:
\jwjwjw{\begin{align*}
    \dot{\mathcal{V}}_\mathcal{N} = &-(x_\mathcal{N}-x_\mathcal{N}^*)^T\begin{bmatrix} \mathcal{R} \otimes I_2 & 0 \\ 0 & \mathcal{G} \otimes I_2\end{bmatrix}(x_\mathcal{N}-x_\mathcal{N}^*) \\&+ (I_{DQ}-I^{\ast}_{DQ})^T(V_{DQ}-V^\ast_{DQ})
\end{align*}}
\jd{From Assumption \ref{ospass} and Definition \ref{passdefn}, we have the following for all buses $j \in N$: there exist open neighbourhoods $W_j$ of $i^{j*}_{DQ}$ and $X_j$ of $x_j^*$ such that $\dot{\mathcal{V}}_{j}(x_{j}) \leq - (i^j_{DQ}-i^{j*}_{DQ})^T(v^j_{DQ}-v^{j*}_{DQ})$} \jd{for all $\il{i^{j}_{DQ}} \in W_j$ and $x_j \in X_j$. This} implies that:
\il{\begin{align*}
    \sum_{j=1}^{|N|}\dot{\mathcal{V}}_{j}(x_{j}) \leq &- (I_{DQ}-I^{\ast}_{DQ})^T(V_{DQ}-V^\ast_{DQ})
\end{align*}}
\jwjwjw{Hence, we have:}
\jwjwjw{\il{\begin{align}
     \dot{\mathcal{V}} \leq -(x_\mathcal{N}-x_\mathcal{N}^*)^T\begin{bmatrix} \mathcal{R} \otimes I_2 & 0 \\ 0 & \mathcal{G} \otimes I_2\end{bmatrix}(x_\mathcal{N}-x_\mathcal{N}^*) \label{lyapfinal}
\end{align}}}
Each $\mathcal{V}_j(x_j)$ has a strict local minimum at $x_j^*$, 
and $\mathcal{V}_\mathcal{N}$ is minimized at \jwjwjw{$(x_\mathcal{N}^\ast)$}. \jd{Thus $\mathcal{V}$ has a strict local minimum at the point $q^* := (\jwjwjw{x_\mathcal{N}^\ast}, x_1^*, x_2^*, ..., x_{|N|}^*)$. We then choose a neighbourhood about $q^*$ on which the following all hold:}}
\begin{enumerate}
    \item $q^*$ is a strict minimum of $\mathcal{V}$. 
    \item $i^{j*}_{DQ} \in W_j$ and $x_j \in X_j$ for all $j \in N$. 
    \item $x_j \in X_j^0$ for all $j \in N$. 
\end{enumerate}
\jd{Thus from \eqref{lyapfinal}, $\mathcal{V}$ is \il{nonincreasing with time}
and has a strict local minimum at $q^*$. Hence,
\il{for sufficiently small $\epsilon>0$ the level set $T=\{q:\mathcal{V}(q)-\mathcal{V}(q^\ast) \leq \epsilon\}$ is
compact and positively invariant with respect to the system \eqref{netwdyn}-\eqref{eq:ss}.}
We then apply LaSalle's Theorem \cite{khalil1991} with \il{function $\mathcal{V}$ on the compact positively invariant set $T$ as the Lyapunov-like function. This states} that all solutions with initial conditions within $T$ converge to the largest invariant set $\Xi$ within $T$ that satisfies $\dot{\mathcal{V}} = 0$.
From \eqref{lyapfinal} $\dot{\mathcal{V}} = 0$ implies \il{that
$x_\mathcal{N} = x_\mathcal{N}^\ast$ and thus $I^N_{DQ} = I^{N\ast}_{DQ}$ \jwjwjw{and $V^{C}_{DQ} = V^{C\ast}_{DQ}$. The former further implies that }}
$I_{DQ} = I^*_{DQ}$ from \eqref{netwdyn}.
Since the inputs to the bus dynamical systems \il{$w_j^* = i^{j*}_{DQ}$} are constant within $\Xi$, \il{we have that each trajectory in $\Xi$ converges to the equilibrium point $q^\ast$ as $t\to\infty$, as follows from the asymptotic stability of $x_j^\ast$ for constant $w_j=w_j^\ast$ (Assumption \ref{ospass}). Since $\mathcal{V}$ has a strict minimum at $q^\ast$ and is constant in $\Xi$ w.r.t. time, we have that $q(t)=q^\ast$ for all $t$ is the only trajectory in $\Xi$.}
Hence, all solutions with initial conditions within $T$ converge to the equilibrium point $q^* = (\jwjwjw{x_\mathcal{N}^\ast}, x_1^*, ..., x_{|N|}^*)$.}

\subsection{Proof of Theorem~\ref{thmlmi}}
\jdwrev{The dynamics considered (the state feedback controller $u = -K\tilde{x} - Mw$ in conjunction with the dynamics~\eqref{eq:convdyn2}) are linear. Hence without loss of generality (and for notational convenience) we can check the passivity property \li{about} the origin, since the passivity property of a linear system is not dependent on its equilibrium point. We will also neglect the bus index $j$ for simplicity, and use the following definitions as in~\eqref{lmi}: $A_c = \tilde{A}-\tilde{B}_uK$, $B_c = \tilde{B}_w - \tilde{B}_uM$, $C_c = \tilde{C} - \tilde{D}_uK$, and $D_c = \tilde{D}_w - \tilde{D}_uM$. Consider \licc{the storage function}
\begin{equation}\label{eq:storage}
    \mathcal{V} = \frac{1}{2}\tilde{x}^TP\tilde{x}
\end{equation}
\jdwrev{in conjunction with the system~\eqref{eq:convdyn2} and $u = -K\tilde{x} - Mw$. To satisfy Assumption~\ref{ospass} \li{with output-strict passivity index $\rho$, we require from \eqref{eq:osp} in \lic{Definition~\ref{passdefn}, }}
$\dot{\mathcal{V}} \leq w^Tz - \rho z^T z $.}
\jdwr{The RHS can be expanded as:}
\licc{\begin{align}
    w^Tz - \rho z^T z =& \frac{1}{2}[w^T(C_c\tilde{x}+D_cw) + (C_c\tilde{x}+D_cw)^Tw] - \nonumber\\&-\rho[(C_c\tilde{x}+D_cw)^T(C_c\tilde{x}+D_cw)] \label{osp}
\end{align}}
\jdwrev{and $\dot{\mathcal{V}}$ is obtained by differentiating $\mathcal{V}$ \jdwr{along the trajectories of~\eqref{eq:convdyn2}, noting that $u = -K\tilde{x} - Mw$:}}
\begin{equation}
\begin{split}
    &\dot{\mathcal{V}} = \frac{1}{2}[\dot{\tilde{x}}^TP\tilde{x}+\tilde{x}^TP\dot{\tilde{x}}]\\
    &~~= \frac{1}{2}[(A_c\tilde{x} + B_cw)^TP\tilde{x}+\tilde{x}^TP(A_c\tilde{x} + B_cw)]\\
    &~~= \frac{1}{2}[\tilde{x}^T(A_c^TP+PA_c)\tilde{x}+\tilde{x}^TPB_c\jwwj{w+w^TB_c}^TP\tilde{x}]
\end{split}\label{vdot}
\end{equation}
\jdwrev{We now recall the condition in Theorem 2 and the associated equation~\eqref{lmi}, and note that if~\eqref{lmi} holds, so does:}
\begin{equation}\label{lmi2}{\footnotesize	
    \begin{bmatrix}
    \tilde{x} \\ w
    \end{bmatrix}^T \begin{bmatrix}	
    A_c^TP+PA_c+2\rho C_c^TC_c & PB_c - C_c^T + 2\rho C_c^TD_c \\	
    B_c^TP - C_c + 2\rho D_c^TC_c & 2\rho D_c^TD_c-D_c^T-D_c	
    \end{bmatrix} \begin{bmatrix}
    \tilde{x} \\ w
    \end{bmatrix}	
    \leq 0}
\end{equation}
\jdwr{From ~\eqref{osp} and~\eqref{vdot}, $\dot{\mathcal{V}} \leq w^Tz - \rho z^T z $ holds if \eqref{lmi2} holds. Hence,
satisfying Theorem~\ref{thmlmi} guarantees that \lic{the passivity condition in} Assumption~\ref{ospass} is satisfied.}}

\lic{
Regarding the asymptotic stability for constant input $w$ we consider the case $w=0$ without loss of generality (since the system is linear).
We have in this case from the analysis above that $\dot{\mathcal{V}} \leq -\rho z^T z $. Asymptotic stability then follows \licc{from Lasalle's Theorem} with $\mathcal{V}$ as the Lyapunov function. In particular, $\dot{\mathcal{V}}=0$ implies $z=0$. From the system dynamics in \eqref{eq:convdyn} we have in this case $i_{iDQ}=0$ and from \eqref{eq:int} $\zeta$ is \licc{constant\footnote{\licc{It should be noted that the constant $\zeta$ is also unique if the feedback gain associated with $\zeta$ is non-zero, a property that should be satisfied by a reasonable controller (otherwise state $\zeta$ can be removed from the controller).}}.} Since the Lyapunov function is radially unbounded an arbitrarily large compact invariant set can be constructed in which the largest invariant set for which $\dot{\mathcal{V}}=0$ is the equilibrium point, thus deducing asymptotic stability. }


\subsection{LMI formulation of the $H_\infty$ and passivity conditions with dynamical controller}
\kl{In this appendix, we present an LMI formulation  for the \li{problem of finding a controller such that we have} output strict passivity of $T_{-{i}_{DQ} \rightarrow v_{DQ}}$ as in~\eqref{eq:osp} together with the $H_\infty$ condition in~\eqref{eq:Hinf_con}.}

\begin{figure}[!t]
	\centering
	\includegraphics[width=0.88\linewidth]{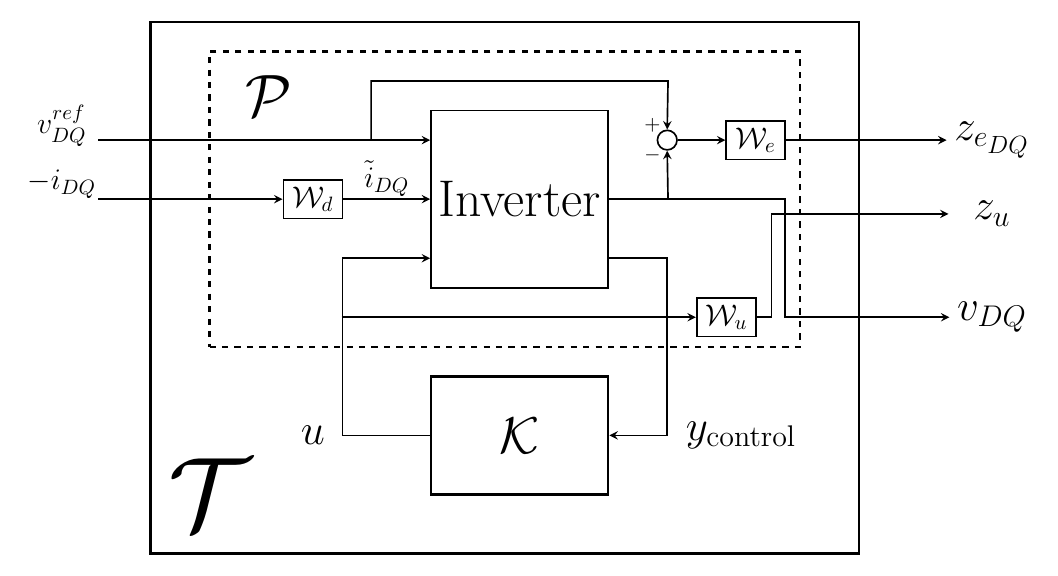}
	\caption{\kkl{\li{Block diagram representing the  $H_\infty$/passivity control synthesis problem. The inputs and outputs of the block $\mathcal{T}$ are associated with the transfer functions on which $H_\infty$/passivity conditions are imposed.}}}
	\label{fig:generalizedplant}
\end{figure}

\kkl{\li{As is commonly followed in 
$H_\infty$ control, to facilitate the control synthesis problem, we consider the block diagram in  
Fig.~\ref{fig:generalizedplant} which contains: the closed loop system $\mathcal{T}$ which has inputs and outputs associated with the transfer functions on which performance specifications are imposed}, the controller $\mathcal{K}$ (given in~\eqref{eq:dynamic_controller}) and the generalized plant $\mathcal{P}$ (containing the inverter itself and the different weighting filters~$W_e$, $W_u$ and~$W_d$).}

\kkl{The different performance specifications (passivity and $H_\infty$) will be expressed using the performance input and output   $w_{\text{perf}}=[(v_{DQ}^{ref})^T~ -{i}_{DQ}^T]^T$ and $z_{\text{perf}}=[	z_{e_{DQ}}^T ~ z_{u}^T~ v_{DQ}^T]^T$. \li{In particular,} 
\begin{itemize}
	\item For passivity,   the input  and output signals are   $w_{\text{pas}}=-{i}_{DQ}$  and   $z_{\text{pas}}=v_{DQ}$.
	\item For $H_\infty$ constraints:  the input  and output signals are   $w_{\infty}=[	(v_{DQ}^{ref})^T~ -{i}_{DQ}^T]^T$  and   $z_{\infty}=[  	z_{e_{DQ}}^T ~ z_{u}^T ]^T$
\end{itemize}}

%

\kkl{To formulate these performance specifications, \li{a state space} representation of the generalized plant $\mathcal{P}$ is required. \jjwr{We} start with the weighting filters.
As explained in \jjwr{Section}~\ref{sec:HinfPass_case_study}, the weighting filters  $W_e$, $W_u$ and $W_d$ are
added to account for loop shaping constraints on $T_{v^{ref}_{DQ} \rightarrow e_{DQ} }$, $T_{v^{ref}_{DQ} \rightarrow u }$, $T_{-i_{DQ} \rightarrow v_{DQ} }$ and $T_{-i_{DQ} \rightarrow u }$. \li{We denote $A_{W_e},B_{W_e},C_{W_e},D_{W_e}$ the matrices in the state space representation \lic{of 
$W_e$ (\jwwj{similarly} for $W_u$ and $W_d$).}} 
The state vectors of  $W_e$, $W_u$ and $W_d$ are \li{denoted as} $x_{W_e}$, $x_{W_u}$ and $x_{W_d}$ respectively. Their inputs are $(v^{ref}_{DQ}-v_{DQ})$, $u$ and $-i_{DQ}$ respectively while their outputs are denoted by $z_{e_{DQ}}$, $z_{u}$ and~$\tilde{i}_{DQ}$ respectively \li{(see Fig.~\ref{fig:generalizedplant})}.}

\kkl{\li{Using}  $\mathcal{A}$, $\mathcal{B}_u$ and $\mathcal{B}_w$ given by~\eqref{eq:ssABCD},  the generalized plant $\mathcal{P}$   has the following representation
\begin{equation}
	\mathcal{P}: ~~~
	\begin{bmatrix}
		\dot{\xi}\\   z_{\text{perf}}\\y_{\text{control}}
	\end{bmatrix}=
	\begin{bmatrix}
		A_\xi &  B_1 & B\\
		C_1 &  D_{11} &  D_{12}\\
		C &  D_{21} &  \bf{0}\\
	\end{bmatrix}
	\begin{bmatrix}
		{\xi}\\   w_{\text{perf}}\\u
	\end{bmatrix}
	\label{eq:generalized_plant}
\end{equation}		where 	 ${\small \xi=[	x^T ~ x_{W_e}^T ~ x_{W_u}^T ~ x_{W_d}^T]^T}$;	  ${\small w_{\text{perf}}=[ 	(v_{DQ}^{ref})^T ~  -{i}_{DQ}^T]^T}$;\\	
${\small z_{\text{perf}}=[  	z_{e_{DQ}}^T ~ z_{u}^T ~ v_{DQ}^T]^T}$;
${\small y_{\text{control}}=[ 	(v^{ref}_{{DQ}})^T ~ \tilde{i}_{DQ}^T ~ i_{iDQ}^T ~ v_{DQ}^T]^T};$
%
\begin{equation}
	{\small \begin{split}
			&A_\xi=\begin{bmatrix}~
				\mathcal{A} & \bf{0} &\bf{0} &C_{W_d}\\ -[\bf{0}~~ B_{W_e}] & A_{W_e}&\bf{0}&\bf{0}\\\bf{0}&\bf{0}&A_{W_u}&\bf{0}\\\bf{0}&\bf{0}&\bf{0}&A_{W_d}
			\end{bmatrix}  \\&B_1=\begin{bmatrix}
				\bf{0}& \mathcal{B}_{w}D_{W_d}\\B_{W_e}&\bf{0}\\\bf{0}&\bf{0}\\\bf{0}&B_{W_d}
			\end{bmatrix}~~~~B=\begin{bmatrix}
				\mathcal{B}_{u}\\\bf{0}\\B_{W_u}\\\bf{0}
			\end{bmatrix}\\&C_1=\begin{bmatrix}
				-[\bf{0} ~~D_{W_e}]& C_{W_e}&\bf{0} &\bf{0}\\\bf{0}&\bf{0}&C_{W_u}&\bf{0}
			\end{bmatrix}\\&  C=\begin{bmatrix}
				I_4 &\bf{0}&\bf{0}&\bf{0}\\\bf{0}&\bf{0}&\bf{0}&\bf{0}
			\end{bmatrix}\\& D_{11}=\begin{bmatrix}
				D_{W_e} & \bf{0}\\\bf{0}&\bf{0}
			\end{bmatrix},
			D_{12}=\begin{bmatrix}\bf{0}\\D_{W_u}
			\end{bmatrix},D_{21}=\begin{bmatrix}\bf{0}\\I_4
			\end{bmatrix}
	\end{split}}
	\label{eq:generalized_plant_matrices}
\end{equation}}

\kl{We are now able to  present an  LMI formulation  for  the
\li{problem of finding a controller $\mathcal{K}$ such that the required passivity and $H_\infty$ conditions are satisfied.}}

\kkl{\begin{thm}\label{thmHinflmi}
		Consider the dynamic feedback controller in~\eqref{eq:dynamic_controller} in conjunction with the dynamics~\eqref{eq:convdyn}.  \\	
		The system~~\eqref{eq:convdyn}
		with input~$w$ and output~$z$ satisfies 1) the passivity condition in Definition~\ref{passdefn}  referred to in Assumption~\ref{ospass} and is output-strictly passive with strictness index $\eta$ and 2) the~$H_\infty$ condition in~\eqref{eq:Hinf_con}
		if  there exist $X$, $Y$, $\hat{A}$, $\hat{B}$, $\hat{C}$, and $\hat{D}$  such that the following LMIs hold
		\begin{equation}
			{\small	\begin{bmatrix} X & I \\I & Y	\end{bmatrix}>0}~~~~~~~~~~~~~~~~~~~~~~~~~~~~~
			\label{eq:QS_condition}
		\end{equation}
		\begin{equation}
			{\small \begin{bmatrix}
					M_{11} & M_{12} & M_{13_\text{pas}} & M_{14_\text{pas}} \\
					M_{12}^T & M_{22} & M_{23_\text{pas}} & M_{24_\text{pas}} \\
					M_{13_\text{pas}}^T & M_{23_\text{pas}}^T & M_{33_\text{pas}} & M_{34_\text{pas}} \\
					M_{14_\text{pas}}^T & M_{24_\text{pas}}^T & M_{34_\text{pas}}^T & -\frac{1}{2}\eta^{-1}I	 		 	
				\end{bmatrix}<0}
			\label{eq:pass_LMI}
		\end{equation}
		\begin{equation}
			{\small	\begin{bmatrix}
					M_{11} & M_{12} & M_{13_\infty} & M_{14_\infty} \\
					M_{12}^T & M_{22} & M_{23_\infty} & M_{24_\infty} \\
					M_{13_\infty}^T & M_{23_\infty}^T & -I & M_{34_\infty} \\
					M_{14_\infty}^T & M_{24_\infty}^T & M_{34_\infty}^T & -I 	 		 	
				\end{bmatrix}<0~~~}
			\label{eq:HInf_LMI}
		\end{equation}			
		with
		\begin{equation*}
			\hspace{-0.35cm}
			{\small\begin{split}	
					&~~M_{11}=A_\xi X+XA_\xi^T+B\hat{C}+(B\hat{C})^T\\
					&~~M_{12}=\hat{A}^T+(A_\xi+B\hat{D}C)\\	
					&~~M_{22}=A_\xi^TY+YA_\xi+\hat{B}C+(\hat{B}C)^T\\
					&~~M_{13_\text{pas}}=(B_{1}R_\text{pas}+B\hat{D}D_{21}R_\text{pas})-(L_\text{pas} C_{1} X+L_\text{pas} D_{12}\hat{C})^T\\
					&~~M_{14_\text{pas}}=(L_\text{pas} C_{1} X+L_\text{pas} D_{12}\hat{C})^T\\	
					&~~M_{23_\text{pas}}=(YB_{1}R_\text{pas}+\hat{B}D_{21}R_\text{pas})-(L_\text{pas} C_{1}+L_\text{pas} D_{12}\hat{D}C)^T\\
					&~~M_{24_\text{pas}}=(L_\text{pas} C_{1}+L_\text{pas} D_{12}\hat{D}C)^T\\
					&~~M_{33_\text{pas}}=-(L_\text{pas} D_{11}R_\text{pas} +L_\text{pas} D_{12}\hat{D}D_{21}R_\text{pas})\\&~~~~~~~~~~~~~~~~-(L_\text{pas} D_{11}R_\text{pas} +L_\text{pas} D_{12}\hat{D}D_{21}R_\text{pas})^T\\
					&~~M_{34_\text{pas}}=(L_\text{pas} D_{11}R_\text{pas} +L_\text{pas} D_{12}\hat{D}D_{21}R_\text{pas})^T\\
					&~~M_{13_\infty}=(B_1R_\infty+B\hat{D}D_{21}R_\infty)\\	
					&~~M_{14_\infty}=(R_\infty C_{1} X+L_\infty D_{12}\hat{C})^T\\	
					&~~M_{23_\infty}=(YB_1R_\infty+\hat{B}D_{21}R_\infty)\\
					&~~M_{24_\infty}=(R_\infty C_{1}+L_\infty D_{12}\hat{D}C)^T\\	
					&~~M_{34_\infty}=(L_\infty D_{11}R_\infty+L_\infty D_{12}\hat{D}D_{21}R_\infty)^T\\
			\end{split}}
		\end{equation*}
		where $A_\xi$, $B$, $C$ $B_1$, $D_{11}$, $D_{12}$ and $D_{21}$ are given by~\eqref{eq:generalized_plant_matrices} and
		\begin{equation}
			\begin{matrix}
				&L_{\text{pas}}&=&
				\begin{bmatrix}  	 0_{2\times 4}&I_2 \end{bmatrix} & R_{\text{pas}}&=&
				\begin{bmatrix}0_{2\times 2}&I_2\end{bmatrix}\\
				&L_{\infty}&=&		\begin{bmatrix}	I_4& 0_{4\times 2}\end{bmatrix} &R_{\infty}&=&I_4 	
			\end{matrix}
			\label{eq:IO_channel_sel}
		\end{equation}	
		The controller $\mathcal{K}$ is \li{constructed} as
		\begin{equation}		
			{\small\begin{split}
					&\li{D_K=\hat{D}, \qquad
						C_K=(\hat{C}-D_KCX)M^{-T}}\li{,}\\
					&B_K=N^{-1}(\hat{B}-YBD_K)\li{,}\\
					&A_K=N^{-1}(\hat{A}-NB_KCX-YBC_KM^T)M^{-T}\\&~~~~~~~~~~~~~-N^{-1}Y(A+BD_KC)XM^{-T}
			\end{split}}
		\end{equation}
		where $M$ and $N$ satisfy $MN^T=I-XY$  which can be solved through a singular value decomposition and Cholesky factorization.
\end{thm}}

\kkl{\subsubsection*{Proof}
	Given the generalized plant dynamics~\eqref{eq:generalized_plant} and the 	
	controller~\eqref{eq:dynamic_controller},  \li{we consider the closed-loop system} 
	\begin{equation*}
		\begin{matrix}
			\mathcal{T}: ~~~
			\begin{bmatrix}
				\dot{x}_{cl}\\z_{\text{perf}}
			\end{bmatrix}=
			\begin{bmatrix}
				A_{c\ell}& B_{c\ell}\\C_{c\ell}& D_{c\ell}
			\end{bmatrix}
			\begin{bmatrix}
				x_{c\ell}\\w_{\text{perf}}
			\end{bmatrix}
		\end{matrix}
		\label{eqt:LPV_closed_loop}
	\end{equation*}
	\li{where} ${x}_{cl}=\begin{bmatrix}
		\xi^T & x_K^T
	\end{bmatrix}^T$ and
	{\small \begin{align*}
			&A_{c\ell}=\begin{bmatrix}
				A_\xi+B D_K C & B C_K \\B_K C & A_K
			\end{bmatrix}
			&B_{c\ell}=\begin{bmatrix}
				B_1 +B D_K  D_{21} \\ B_K  D_{21} 	
			\end{bmatrix} \\
			&C_{c\ell}=\begin{bmatrix}
				C_1+  D_{12}D_KC&  D_{12}C_K	
			\end{bmatrix}
			&D_{c\ell}=  D_{11}+  D_{12}D_K  D_{21}
	\end{align*}} To formulate the passivity condition and the $H_\infty$ condition, we extract the corresponding closed-loop systems \li{on which these conditions are imposed}
	by selecting the appropriate input-output channels \li{of $\mathcal{T}$} using  $L_{\text{pas}}$, $R_{\text{pas}}$,  $L_{\infty}$ and $R_{\infty}$ (given by~\eqref{eq:IO_channel_sel}) as \li{follows\footnote{\li{$A_{\text{pas}}, B_{\text{pas}}, C_{\text{pas}}, D_{\text{pas}}$ denote the matrices in the state space representation of system $\mathcal{T}_{\text{pas}}$ (similarly for $\mathcal{T}_{\infty}$).}}}
	\begin{equation}
		\mathcal{T}_{\text{pas}}:{\small \begin{bmatrix}
				A_{\text{pas}}& B_{\text{pas}} \\  C_{\text{pas}}&   D_{\text{pas}}
		\end{bmatrix}}= 		{\small \begin{bmatrix}
				A_{c\ell}& B_{c\ell} R_{\text{pas}} \\ L_{\text{pas}} C_{c\ell}&  L_{\text{pas}} D_{c\ell} R_{\text{ps}}
		\end{bmatrix}}~~~~~~~~~~     \label{eq:closedLopp_Pass_matrices}
	\end{equation}	
	\begin{equation}\hspace{-0.35cm}
		~~~~~~~  	\mathcal{T}_{\infty}:{\small \begin{bmatrix}
				A_{\infty}& B_{\infty} \\  C_{\infty}&   D_{\infty}
		\end{bmatrix}}={\small \begin{bmatrix}
				A_{c\ell}& B_{c\ell} R_{\infty}\\ L_{\infty} C_{c\ell}&  L_{\infty} D_{c\ell} R_{\infty}
		\end{bmatrix}}
		\label{eq:closedLopp_Hinf_matrices}
	\end{equation}
	Using standard results from LQ control theory~\cite{Willems1971}, the  output-strict passivity of $\mathcal{T}_{\text{pas}}$ with passivity index~$\eta$ is satisfied if and only if  there exists $W_\text{pas}=W_\text{pas}^T>0$ such that
	\begin{equation}
		{\small	\begin{split}
				& \begin{bmatrix}		A_{\text{pas}}^T W_{\text{pas}} + W_{\text{pas}}A_{\text{pas}} & W_{\text{pas}} B_\text{pas}-C_\text{pas}^T & C_\text{pas}^T \\
					B_\text{pas}^T W_{\text{pas}}-C_\text{pas} & -D_\text{pas}^T-D_\text{pas}  & D_\text{pas}^T	\\
					C_\text{pas}& D_\text{pas}& -\frac{\eta^{-1}}{2}I  \end{bmatrix} <0
		\end{split}}
		\label{eq:PassBMI}
	\end{equation} 	
	In the same way, the condition		
	$\left\|\mathcal{T}_{\infty} \right\|_\infty \leq 1$  is satisfied if and only if there exists $W_{\infty}=W_{\infty}^T>0$ such that
	\begin{equation}
		{\small \begin{split}
				&\begin{bmatrix}		A_{\infty}^T W_{\infty} + W_{\infty} A_{\infty} & W_{\infty} B_{\infty} & C_{\infty}^T \\
					B_{\infty}^T W_{\infty} & -I&	D_{\infty}^T \\
					C_{\infty} & D_{\infty} & -I
				\end{bmatrix} <0~
		\end{split}}
		\label{eq:HinfBMI}
	\end{equation}
	We use \li{the results in~\cite{scherer1997}} and~\cite[Chap 4]{ScW:05} to formulate  LMIs through which $A_K$, $B_K$, $C_K$ and $D_K$  can be found such that conditions~\eqref{eq:PassBMI}-\eqref{eq:HinfBMI} are satisfied.\\
	\li{In particular, let\footnote{\lif{This will make conditions~\eqref{eq:PassBMI}-\eqref{eq:HinfBMI} only sufficient.} \li{As noted in \cite{scherer1997} this allows to achieve an LMI formulation for an otherwise non-convex problem.}} $W_{\infty}=W_\text{pas}=W$ and} 
	$${\small
		W=\begin{bmatrix} 	Y & N \\ N^T & U \end{bmatrix},~~W^{-1}=\begin{bmatrix} 	X & M \\ M^T & V \end{bmatrix},~~\Phi=\begin{bmatrix}
			X & I\\M^T&\bf{0}
	\end{bmatrix}}$$
	where $X$, $Y$, $V$ and $U$ are symmetric matrices. \\
	Due to  $WW^{-1}=I$, we have  $YX+NM^T=I$, $N^TX+UM^T=\bf{0}$, $YM+NV=\bf{0}$ and $UV-XY=\bf{0}$. Therefore,  pre and post-multiplying of $W>0$ by $\Phi$ \jjwr{results} in condition~\eqref{eq:QS_condition}.\\
	We define now the \li{following change of 
		variables} 
	\begin{equation}
		\hspace{-0.25cm}
		{\small	\begin{split}
				&\hat{A}=NA_KM^T+NB_KCX+YBC_KM^T+Y(A+BD_KC)X\\
				&\li{\hat{B}=NB_K+YBD_K, \quad
					\hat{C}=C_KM^T+D_KCX, \quad
					\hat{D}=D_K}
		\end{split}}
		\label{eq:CVK}
	\end{equation} Conditions~\eqref{eq:pass_LMI}-\eqref{eq:HInf_LMI} are obtained by pre and postmultiplying  conditions~\eqref{eq:PassBMI}-\eqref{eq:HinfBMI} (with ${\small W_{\infty}=W_\text{pas}=W}$) by ${\small \textbf{bdiag}(\Phi,I,I)}$ and using~\eqref{eq:CVK} to rewrite the obtained conditions\jjwr{, \li{thus} completing the proof}. }

\end{document}